\documentclass[preprint,3p,12pt]{elsarticle}
\usepackage{mathrsfs}
\usepackage{amsmath}
\usepackage{stmaryrd}
\usepackage{bbding}
\usepackage{dcolumn}
\usepackage{graphicx}
\usepackage{amsfonts}
\usepackage{amssymb}
\usepackage{psfrag}
\usepackage{wrapfig}
\usepackage{makeidx}
\usepackage{bm}
\usepackage{epsf}
\usepackage{epsfig}
\usepackage{setspace}
\usepackage{graphicx}
\usepackage{epstopdf}
\usepackage{psfrag}
\usepackage{subcaption}
\usepackage{color}
\usepackage{comment}
\usepackage{hyperref}

\begin{document}
	
	\title{Comparison of the high-order Runge-Kutta discontinuous Galerkin method and gas-kinetic scheme for inviscid compressible flow simulations}
	
	\author[XJTU]{Yixiao Wang}
	\ead{2213411837@stu.xjtu.edu.cn}
	
	\author[XJTU]{Xing Ji\corref{cor1}}
	\ead{jixing@xjtu.edu.cn}
	
	\author[XJTU]{Gang Chen}
	\ead{aachengang@xjtu.edu.cn}
	
	\author[HKUST1,HKUST2,HKUST3]{Kun Xu}
	\ead{makxu@ust.hk}

	\address[XJTU]{Shaanxi Key Laboratory of Environment and Control for Flight Vehicle, Xi'an Jiaotong University, Xi'an, China}
	\address[HKUST1]{Department of Mathematics, Hong Kong University of Science and Technology, Clear Water Bay, Kowloon, Hong Kong}
	\address[HKUST2]{Department of Mechanical and Aerospace Engineering, Hong Kong University of Science and Technology, Clear Water Bay, Kowloon, Hong Kong}
	\address[HKUST3]{Shenzhen Research Institute, Hong Kong University of Science and Technology, Shenzhen, China}
	\cortext[cor1]{Corresponding author}
	
	\begin{abstract}
		
The Runge-Kutta Discontinuous Galerkin (RKDG) method is a high-order technique for addressing hyperbolic conservation laws, which has been refined over recent decades and is effective in handling shock discontinuities. Despite its advancements, the RKDG method faces challenges, such as stringent constraints on the explicit time-step size and reduced robustness when dealing with strong discontinuities. On the other hand, the Gas-Kinetic Scheme (GKS) based on a high-order gas evolution model also delivers significant accuracy and stability in solving hyperbolic conservation laws through refined spatial and temporal discretizations. Unlike RKDG, GKS allows for more flexible CFL number constraints and features an advanced flow evolution mechanism at cell interfaces. Additionally, GKS’s compact spatial reconstruction enhances the accuracy of the method and its ability to capture stable strong discontinuities effectively. In this study, we conduct a thorough examination of the RKDG method using various numerical fluxes and the GKS method employing both compact and non-compact spatial reconstructions. Both methods are applied under the framework of explicit time discretization and are tested solely in inviscid scenarios. We will present numerous numerical tests and provide a comparative analysis of the outcomes derived from these two computational approaches.

	\end{abstract}
	
	
	\maketitle
	
	\section{Introduction}

In this paper, we address the solutions of one- and two-dimensional nonlinear hyperbolic conservation laws using the Runge-Kutta Discontinuous Galerkin (RKDG) method and the Gas-Kinetic Scheme (GKS), evaluating their performance across various standard test cases. The discontinuous Galerkin method, first proposed by Reed and Hill for neutron transport \cite{reed1973triangular}, has seen substantial developments by Cockburn et al., establishing a robust framework for solving nonlinear, time-dependent hyperbolic conservation laws with explicit, strong stability-preserving Runge-Kutta time discretization methods \cite{shu1988efficient,cockburn1989tvb,cockburn1990runge,cockburn1998runge}. The RKDG method is particularly adept at capturing weak discontinuities and offers considerable flexibility, ease of parallelization, and scalability to high-order cases. However, it may produce significant oscillations when dealing with strong discontinuities, particularly in high-order applications.
To mitigate these oscillations, several nonlinear limiters have been introduced. Classical minmod-type Total Variation Bounded (TVB) limiters help restrict the solution's slope \cite{cockburn1989tvb}, yet they may compromise accuracy in smoothly varying regions by inaccurately identifying trouble cells. Recently, the oscillation-free Discontinuous Galerkin (OFDG) methods have been developed \cite{liu2022essentially,lu2021oscillation} , incorporating a damping term proportional to the discontinuity's intensity to better control spurious oscillations. These methods, however, often require tighter CFL restrictions or the adoption of exponential Runge-Kutta time discretization strategies.
Additionally, limiters based on Weighted Essentially Non-Oscillatory (WENO) and High-Order Weighted Essentially Non-Oscillatory (HWENO) reconstructions have successfully merged the high accuracy of classical WENO methods with the compactness of the DG method \cite{jiang1996efficient,friedrich1998weighted,hu1999weighted,qiu2005runge,qiu2005hermite,zhu2008runge,zhu2016runge}, albeit at the cost of increased computational demands. Among these, the multi-resolution WENO limiter stands out for its compact spatial stencil and straightforward implementation \cite{zhu2020high,zhu2018new}, making it a promising choice for inclusion in our study.

Over recent years, the Gas-Kinetic Scheme (GKS) has undergone systematic development and has proven effective in solving Euler, Navier-Stokes (N-S) flows, and beyond \cite{xu2014direct,xu2001gas,xu2010unified}. GKS is grounded in gas kinetic theory at the mesoscopic scale, utilizing kinetic equations such as the Bhatnagar-Gross-Krook (BGK) model to describe mesoscopic gas particle evolution \cite{bhatnagar1954model}. The macroscopic gas dynamic equations, including the N-S and Euler equations, are derived from the BGK model using the Chapman-Enskog expansion  \cite{xu2014direct,chapman1990mathematical}, which helps in determining the time evolution of gas distribution function.
A distinctive feature of GKS is its approach in the flux evaluation at cell interfaces. Unlike traditional numerical fluxes, GKS provides both the flux value and its time derivative, along with time-accurate flow variables at the cell interface. This method has demonstrated high accuracy in smooth flow regions while introducing necessary numerical dissipation to maintain robustness in discontinuous areas, as highlighted in previous studies \cite{ji2018family,ji2018compact}. The introduction of a two-stage fourth-order time discretization framework \cite{pan2016efficient}, along with compact spatial discretization techniques, has further enhanced the efficiency of the Compact Gas-Kinetic Scheme (CGKS). These advancements make CGKS particularly adept at handling complex shock interactions, positioning it as a highly effective tool in computational fluid dynamics simulations.

The primary differences between the RKDG method and the GKS are outlined as follows:

\noindent 1. CFL Number Flexibility: GKS operates with fewer restrictions on the Courant-Friedrichs-Lewy (CFL) number compared to RKDG. However, reaching the accuracy level of RKDG or other high-order schemes can be challenging for GKS, mainly due to the constraints imposed by its finite volume framework.

\noindent 2. Identification of Troubled Cells and Limiting: RKDG and GKS use markedly different approaches for detecting troubled cells and applying limits to the numerical solution. These differences stem from their unique ways of handling the evolution of degrees of freedom within each cell.

\noindent 3. Mathematical Formulation and Physical Foundation: RKDG is known for its straightforward mathematical formulation, though its performance can be problematic in regions with strong discontinuities, where robustness tends to deteriorate quickly. In contrast, GKS, while more complex in formulation and less straightforward to extend to very high orders, is rooted in a physical model of gas dynamics \cite{xu-liu}, offering enhanced flexibility for managing complex flow scenarios.

\noindent 4. Degrees of Freedom and Computational Demands: RKDG allows for high accuracy with limited mesh refinement due to ample degrees of freedom. However, increasing the order of the scheme typically reduces the permissible time step size, demanding greater memory bandwidth. On the other hand, the CGKS may have fewer degrees of freedom but achieves high-order accuracy through spatial reconstruction. Its spatio-temporal coupling approach relaxes CFL conditions for time advancement, although the computation of time-accurate gas-kinetic flux terms can reduce its overall efficiency.

In this paper, we systematically compare these two advanced high-order methods by solving 1-D and 2-D Euler equations for the first time. We note that recent developments, such as the Reconstructed Discontinuous Galerkin (RDG) methods and Arbitrary DERivative Discontinuous Galerkin (ADER-DG) methods \cite{luo2010reconstructed,luo2012reconstructed,kaser2006arbitrary,dumbser2006arbitrary}, have addressed memory overhead in the DG framework. While we continue to refine the efficiency of GKS, this study focuses exclusively on the classical implementations of both methods. To ensure fairness in our comparisons, both methods have been implemented by us using the same data structures.

	The structure of this paper is outlined as follows.
Section 2 provides an overview of the RKDG method, including a description of the multi-resolution WENO-type limiters applied in both one-dimensional and two-dimensional scenarios. Section 3 introduces the GKS and discusses the techniques used for compact spatial reconstruction within this framework.
Section 4 presents comparative analyses of the two methods, supported by results from numerical tests.
The final section offers concluding remarks and summarizes the key findings of the study.

	\section{A brief review of the RKDG method}
	To solve the hyperbolic conservation laws, the first step is to give a partition of the computational domain, including cells $ I_j=\left[x_{j-\frac{1}{2}}, x_{j+\frac{1}{2}}\right], \quad j=1, \cdots, N$, where the cell center is denoted by $x_j=\frac{1}{2}\left(x_{j-\frac{1}{2}}+x_{j+\frac{1}{2}}\right)$, and the mesh size by $\Delta x_j=x_{j+\frac{1}{2}}-x_{j-\frac{1}{2}}$. The solution, as well as the test function space for the Galerkin method, is given by $V_h^k=\left\{p:\left.p\right|_{I_j} \in P^k\left(I_i\right)\right\}$, where $P^k\left(I_i\right)$ is the space of polynomials of degree $\leqslant k$ on the cell $I_i$. In this paper, A local orthogonal basis over $I_i$ $\left\{v_l^{(j)}(x), l=0,1, \ldots, k\right\}$ is used to avoid calculating mass term for better computational efficiency:

	\begin{equation}
		\begin{aligned}
			& v_0^{(j)}(x)=1, \\
			& v_1^{(j)}(x)=\sqrt{12}\left(\frac{x-x_j}{\Delta x_j}\right), \\
			& v_2^{(j)}(x)=\sqrt{180}\left(\left(\frac{x-x_j}{\Delta x_j}\right)^2-\frac{1}{12}\right), \\
			& v_3^{(j)}(x)=\sqrt{2800}\left(\left(\frac{x-x_j}{\Delta x_j}\right)^3-\frac{3}{20}\left(\frac{x-x_j}{\Delta x_j}\right)\right), \\
			& v_4^{(j)}(x)=\sqrt{44100}\left(\left(\frac{x-x_j}{\Delta x_j}\right)^4-\frac{3}{14}\left(\frac{x-x_j}{\Delta x_j}\right)^2+\frac{3}{560}\right)\\
			&\qquad \qquad \qquad  \qquad \qquad  \ldots
		\end{aligned}
	\end{equation}
	
	The one-dimensional solution $u_h(x, t) \in V_h^k$ can be written as:
	\begin{equation}
		u_h(x, t)=\sum_{l=0}^k u_j^{(l)}(t) v_l^{(j)}(x), \quad x \in I_j,
	\end{equation}
	
	and the semi-discrete scheme is introduced to find a unique function $u_h(\cdot, t) \in V_h^k$
	such that for all discrete cells,
	\begin{equation}
		\int_{I_j}\left(u_h\right)_t v d x-\int_{I_j} f\left(u_h\right) v_x d x+\hat{f}_{j+\frac{1}{2}} v\left(x_{j+\frac{1}{2}}^{-}\right)-\hat{f}_{j-\frac{1}{2}} v\left(x_{j-\frac{1}{2}}^{+}\right)=0
	\end{equation}
	holds for all test functions  $v \in V_h^k$.
	Same notations for numerical flux as in \cite{zhu2016runge} are adopted.
	
	For the two-dimensional case, the uniform mesh partition in the computational domain is adopted, consisting of $I_{ij} = \left[x_{i-\frac{1}{2}}, x_{i+\frac{1}{2}}\right] \times \left[y_{j-\frac{1}{2}}, y_{j+\frac{1}{2}}\right], \quad i=1,2,\ldots, N_x$ and $j=1,2,\ldots, N_y$. with the mesh size $\Delta x_i=x_{i+\frac{1}{2}}-x_{i-\frac{1}{2}},\Delta y_j= y_{j+\frac{1}{2}}-y_{j-\frac{1}{2}}$.
	The solution as well as test function space is defined by $S_h^k=\left\{v(x, y):\left.v(x, y)\right|_{I_{i, j}} \in \mathbb{P}^k\left(I_{i, j}\right)\right\}$ as the piecewise polynomials space of degree at most k defined on $I_{ij}$. A local orthonormal basis over $I_{ij}$, $\left\{v_l^{(i, j)}(x, y), l=0,1, \ldots, K ; K=\frac{(k+1)(k+2)}{2}-1\right\}$ is adopted:
	\begin{equation}
		\begin{aligned}
			& v_0^{(i, j)}(x, y)=1, \\
			& v_1^{(i, j)}(x, y)=v_1^{(i)}(x),\\
			& v_2^{(i, j)}(x, y)=v_1^{(j)}(y), \\
			& v_3^{(i, j)}(x, y)=v_2^{(i)}(x), \\
			& v_4^{(i, j)}(x, y)=v_1^{(i)}(x) v_1^{(j)}(y), \\
			&v_5^{(i, j)}(x, y)=v_2^{(j)}(y), \\
			& v_6^{(i, j)}(x, y)=v_3^{(i)}(x), \\
			& v_7^{(i, j)}(x, y)=v_2^{(i)}(x) v_1^{(j)}(y), \\
			& v_8^{(i, j)}(x, y)=v_1^{(i)}(x) v_2^{(j)}(y), \\
			& v_9^{(i, j)}(x, y)=v_3^{(j)}(y), \\
			& v_{10}^{(i, j)}(x, y)=v_4^{(i)}(x), \\
			& v_{11}^{(i, j)}(x, y)=v_3^{(i)}(x) v_1^{(j)}(y), \\
			& v_{12}^{(i, j)}(x, y)=v_2^{(i)}(x) v_2^{(j)}(y),\\
			& v_{13}^{(i, j)}(x, y)=v_1^{(i)}(x) v_3^{(j)}(y), \\
			& v_{14}^{(i, j)}(x, y)=v_4^{(j)}(y)\\
			&\qquad \qquad   \ldots
		\end{aligned}
	\end{equation}

	The two-dimensional solution $u_h(x,y,t)$ can be written as:
	\begin{equation}
		u_h(x, y, t)=\sum_{l=0}^K u_{i, j}^{(l)}(t) v_l^{(i, j)}(x, y), \quad(x, y) \in I_{i, j},
	\end{equation}
	and a two-dimensional semi-discretization scheme is adopted to obtain a specific function $u_h \in S_h^k$ such that for all discrete cells,
	\begin{equation}
		\begin{aligned}
			& \int_{I_{i, j}}\left(u_h\right)_t v d x d y \\
			= & \int_{I_{i, j}} f\left(u_h\right) v_x d x d y-\int_{I_j} \hat{f}_{i+\frac{1}{2}}(y) v\left(x_{i+\frac{1}{2}}^{-}, y\right) d y+\int_{I_j} \hat{f}_{i-\frac{1}{2}}(y) v\left(x_{i-\frac{1}{2}}^{+}, y\right) d y \\
			& +\int_{I_{i, j}} g\left(u_h\right) v_y d x d y-\int_{I_i} \hat{g}_{j+\frac{1}{2}}(x) v\left(x, y_{j+\frac{1}{2}}^{-}\right) d x+\int_{I_i} \hat{g}_{i-\frac{1}{2}}(x) v\left(x, y_{j-\frac{1}{2}}^{+}\right) d x \\
			& :=L(u_h),
		\end{aligned}
	\end{equation}
	holds for all the test function $v \in S_h^k$, where the notations for numerical flux are referred to \cite{zhu2016runge}. The Gaussian quadrature is used to calculate the integral in practice.
	As for the time discretization, a third-order TVD Runge-Kutta method \cite{shu1988efficient} is adopted:
	\begin{equation}\label{RK3}
		\left\{\begin{array}{l}
			u^{(1)}=u^n+\Delta t L\left(u^n\right), \\
			u^{(2)}=\frac{3}{4} u^n+\frac{1}{4} u^{(1)}+\frac{1}{4} \Delta t L\left(u^{(1)}\right), \\
			u^{n+1}=\frac{1}{3} u^n+\frac{2}{3} u^{(2)}+\frac{2}{3} \Delta t L\left(u^{(2)}\right),
		\end{array}\right.
	\end{equation}
	Thus, a fully discrete scheme in time and space is achieved.

	\subsection{Implementation of limiters}
	
	In this subsection, the implementation of limiters for RKDG is introduced, taking the multi-resolution WENO-type limiters for illustration.
	It is crucial to indicate the trouble cells appropriately, for the indicator would degrade the accuracy in smooth regions if it indicates one mistakenly, or it may cause the program's instability if it omits one that ought to be modified.
	In comparing different trouble indicators \cite{qiu2005comparison}, the KXRCF shock detection technique developed in \cite{krivodonova2004shock} is recommended. It divides the boundary of the target cell $I_{ij}$ into two parts: $\partial I_{ij}^{-}$ and $\partial I_{ij}^{+}$ where, respectively, the flow is into $(\vec{v} \cdot \vec{n}<0)$ and out of $(\vec{v} \cdot \vec{n}>0)$.
	The discontinuity detector is defined as
	\begin{equation}
		\frac{\left|\int_{\partial I_{i, j}^{-}}\left(\left.u_h(x, y, t)\right|_{I_{i, j}}-\left.u_h(x, y, t)\right|_{I_l}\right) d s\right|}{\left.h_{i, j}^R\left|\partial I_{i, j}^{-}\right| \cdot||| \widehat{u_h}(x,y,t)\right|_{\partial I_{i, j}}|||} \geq C_k,
	\end{equation}
	where $C_k$ is a constant set to be 1 in our work, here $I_l$ represents the adjacent cells sharing the same boundary $\partial I_{ij}^{-}$ with the target cell.  $||| \widehat{u_h}(x,y,t) _{\partial I_{i, j}}|||$ is defined of the minimum value of $|u_h(x,y,t)|$ in Gauss integral points along $\partial I_{i, j}$. $h_{ij} $ is the radius of the circumscribed circle in the target cell, the power of which, namely R, is recommended to be taken as 1 for k=1 and 1.5 for k$>$1 case in \cite{fu2017new}.
	
	The following gives a brief introduction of the multi-resolution WENO-type limiter.  This method fully utilizes the higher-order information in the solutions of the troubled cell and takes advantage of the orthonormal basis to simplify the computation cost. It adopts a WENO-Z-type scheme and follows the ideas of the CWENO schemes to construct smooth indicators and weights for cells. Same notations in \cite{zhu2020high} are taken in the following, which can be referred to for more details.
	The key procedure of the algorithm is to reconstruct the modified polynomial $\left.u_h^{\text {new }}\right|_{I_{i, j}}=p^{\text {new }}(x, y)$ on the target cell by
	\begin{equation}\label{evolution}
		p^{\text {new }}(x, y)=\sum_{\ell=\ell_2-1}^{\ell_2} \omega_{\ell, \ell_2} p_{\ell, \ell_2}(x, y), \quad \ell_2=1,2,3,4,
	\end{equation}
	for different-order schemes, the provisional polynomials associated with the original solutions and the weights $\omega_{\ell, \ell_2}$ associated with the WENO-Z-type coefficients are taken below.
	
	Initially, we make projections for solutions of the trouble cells and their four adjacent cells. A series of new polynomials are denoted by
	$q^{(\ell)}_{i,j}(x,y), q^{(\ell)}_{i-1,j}(x,y), q^{(\ell)}_{i+1.j}(x,y), q^{(\ell)}_{i,j+1}(x,y)$ and $q^{(\ell)}_{i,j-1}$ respectively, where $\ell $ takes $0,1,\cdots,k$.  Owing to the application of the specified basis, these can be easily obtained by the solutions contained in cells. For simplicity, the subscripts of $q^{(\ell)}_{i,j}(x,y)$ are omitted.
	
	Secondly, a series of limited polynomials $p_{\ell, \ell}(x, y), \ell=1, \ldots, k$ are obtained, together with $p_{\ell, \ell+1}(x, y), \ell=1, \ldots, k-1$ through
	\begin{equation}
		p_{\ell, \ell}(x, y)=\frac{1}{\gamma_{\ell, \ell}} q^{(\ell)}(x, y)-\frac{\gamma_{\ell-1, \ell}}{\gamma_{\ell, \ell}} p_{\ell-1, \ell}(x, y), \ell=1, \ldots, k
	\end{equation}
	and
	\begin{equation}
		p_{\ell, \ell+1}(x, y)=\omega_{\ell, \ell} p_{\ell, \ell}(x, y)+\omega_{\ell-1, \ell} p_{\ell-1, \ell}(x, y), \ell=1, \ldots, k-1
	\end{equation}
	Here, $\gamma_{\ell-1, \ell}+\gamma_{\ell, \ell}=1$, in which $\gamma_{\ell-1, \ell}$ and $\gamma_{\ell, \ell}$ are determined in advance. It should be pointed out that these two coefficients are vital for the results of the simulations, and they are chosen carefully case by case in our work.
	$\omega_{\ell-1, \ell}+\omega_{\ell, \ell}=1$ yields for consistency at the same time.
	
	Next is to calculate the smoothness indicators for the corresponding polynomials.
	Due to the compactness of this algorithm, only the information of the four adjacent cells is used to construct $\beta_{0,1}:$
	\begin{equation}
		\beta_{0,1}=\min \left(\varsigma_{i, j-1}, \varsigma_{i, j+1}, \varsigma_{i-1, j}, \varsigma_{i+1, j}\right) .
	\end{equation}
	Where, $\varsigma_{m,n}$ is defined as:
	\begin{equation}
		\varsigma_{m,n}=\int_{I_{i, j}}\left(\frac{\partial}{\partial x} q_{m,n}^{(\ell)}(x, y)\right)^2+\left(\frac{\partial}{\partial y} q_{m, n}^{(\ell)} (x, y)\right)^2 d x d y
	\end{equation}
	and $m,n$ refer to iterating over the four adjacent cells.
	Then, the smoothness indicators for high-order polynomials are constructed:
	\begin{equation}
		\beta_{\ell, \ell_2}=\sum_{|\alpha|=1}^\kappa \int_{I_{i, j}}\left(\Delta x_i \Delta y_j\right)^{|\alpha|-1}\left(\frac{\partial^{|\alpha|}}{\partial x^{\alpha_1} \partial y^{\alpha_2}} p_{\ell, \ell_2}(x, y)\right)^2 d x d y, \quad \ell=\ell_2-1, \ell_2 ; \ell_2=1,2,3,4
	\end{equation}
	where $\kappa=\ell, \alpha=\left(\alpha_1, \alpha_2\right)$, and $|\alpha|=\alpha_1+\alpha_2$, respectively.
	Then, the WENO-Z recipe \cite{borges2008improved,castro2011high} is adapted to get the relative difference between the smoothness indicators:
	\begin{equation}
		\tau_{\ell_2}=\left(\beta_{\ell_2, \ell_2}-\beta_{\ell_2-1, \ell_2}\right)^2, \quad \ell_2=1,2,3,4,
	\end{equation}
	and nonlinear weights can be derived from
	\begin{equation}
		\omega_{\ell_1, \ell_2}=\frac{\bar{\omega}_{\ell_1, \ell_2}}{\sum_{\ell=1}^{\ell_2} \bar{\omega}_{\ell, \ell_2}}, \quad \bar{\omega}_{\ell_1, \ell_2}=\gamma_{\ell_1, \ell_2}\left(1+\frac{\tau_{\ell_2}}{\varepsilon+\beta_{\ell_1, \ell_2}}\right), \ell_1=\ell_2-1, \ell_2 ; \ell_2=1,2,3,4 .
	\end{equation}
	Here, $\varepsilon$ is take as $10^{-10}$ in our work.
	Finally, substitute the nonlinear weights for \eqref{evolution} to obtain the reconstructed polynomial $\left.u_h^{\text {new }}\right|_{I_{i, j}}=p^{\text {new }}(x, y)$.
	
	For two-dimensional Euler system equations
	\begin{equation}
		\frac{\partial}{\partial t}\left(\begin{array}{c}
			\rho \\
			\rho \mu \\
			\rho \nu \\
			E
		\end{array}\right)+\frac{\partial}{\partial x}\left(\begin{array}{c}
			\rho \mu \\
			\rho \mu^2+p \\
			\rho \mu \nu \\
			\mu(E+p)
		\end{array}\right)+\frac{\partial}{\partial y}\left(\begin{array}{c}
			\rho \nu \\
			\rho \mu \nu \\
			\rho \nu^2+p \\
			\nu(E+p)
		\end{array}\right)=0,
	\end{equation}a characteristic-wise \cite{zhu2016runge} procedure is applied to the limitations.
	
	When it comes to the time discretization,  the modifications above are done to obtain the $u_h^{\text{new}}$ in each sub-time-step in Eq.\eqref{RK3}, which substitutes the original solution in trouble cells for time evolution.

	\section{A brief review of the gas-kinetic method}
	In this section, the gas-kinetic evolution model is introduced first. Then, the implementation of the spatial reconstruction is briefly reviewed.
	
	\subsection{Gas-kinetic evolution model}

	The two-dimensional gas-kinetic BGK equation \cite{bhatnagar1954model} can be
	written as
	\begin{equation}\label{bgk}
		f_t+\textbf{u}\cdot\nabla f=\frac{g-f}{\tau},
	\end{equation}
	where $f$ is the gas distribution function, $g$ is the corresponding equilibrium state.
	
	The equilibrium state is a Maxwellian distribution
	\begin{equation*}
		\begin{split}
			g=\rho(\frac{\lambda}{\pi})^{\frac{K+2}{2}}e^{\lambda((u-U)^2+(v-V)^2+\xi^2)},
		\end{split}
	\end{equation*}
	where $\lambda =m/2kT $, and $m, k, T$ represents the molecular mass, the Boltzmann constant, and temperature, $K$ is the number of internal degrees of freedom, i.e. $K=(4-2\gamma)/(\gamma-1)$ for two-dimensional flows, and $\gamma$ is the specific heat ratio.
	The collision
	term satisfies the following compatibility condition
	\begin{equation}\label{compatibility1}
		\int \frac{g-f}{\tau}{\bm{\psi}} \text{d}\Xi=0,
	\end{equation}
	where $\bm{\psi}=(1,u,v,\displaystyle \frac{1}{2}(u^2+v^2+\xi^2))$,
	$\text{d}\Xi=\text{d}u\text{d}v\text{d}\xi_1...\text{d}\xi_{K}$, the internal variable $\xi^2=\xi^2_1+\xi^2_2+...+\xi^2_K$. The connections between macroscopic mass $\rho$, momentum ($\rho U, \rho V$), and energy $\rho E$ with the distribution function $f$ are
	\begin{equation}\label{f-to-convar}
		\left(
		\begin{array}{c}
			\rho\\
			\rho U\\
			\rho V\\
			\rho E\\
		\end{array}
		\right)
		=\int \bm{\psi} f \text{d}\Xi.
	\end{equation}
	Similar kinetic models for quantum gases with equilibrium Bose-Einstein or Fermi-Dirac distributions have been constructed and
	studied as well \cite{yang2007high}.
	\color{black}{}
	
	Based on the Chapman-Enskog expansion for the BGK equation
	\cite{xu2014direct}, the gas distribution function in the continuum
	regime can be expanded as
	\begin{align*}
		f=g-\tau D_{\textbf{u}}g+\tau D_{\textbf{u}}(\tau
		D_{\textbf{u}})g-\tau D_{\textbf{u}}[\tau D_{\textbf{u}}(\tau
		D_{\textbf{u}})g]+...,
	\end{align*}
	where $D_{\textbf{u}}={\partial}/{\partial t}+\textbf{u}\cdot
	\nabla$. By truncating on different orders of $\tau$, the
	corresponding macroscopic equations can be derived. In this paper, only the Euler equations are discussed. Thus, the zeroth-order truncation is taken, i.e. $f=g$.
	
	Taking moments of the BGK equation Eq.\eqref{bgk} and integrating with
	respect to space, the semi-discrete finite volume scheme can be
	written as
	\begin{align*}
		\frac{\text{d}\bar{\textbf{W}}_{ij}}{\text{d}t}=-\frac{1}{\Delta x}
		(\textbf{F}_{i+1/2,j}(t)-\textbf{F}_{i-1/2,j}(t))-\frac{1}{\Delta y}
		(\textbf{G}_{i,j+1/2}(t)-\textbf{G}_{i,j-1/2}(t)),
	\end{align*}
	where $\bar{\textbf{W}}_{ij}$ is the cell averaged value of conservative variables,
	$\textbf{F}_{i \pm 1/2,j}(t)$ and $\textbf{G}_{i,j \pm 1/2}(t)$ are the time-dependent
	numerical fluxes at cell interfaces in $x$ and $y$ directions.
	The Gaussian quadrature is used to achieve the accuracy in space, such that
	\begin{align}\label{gauss-flux}
		\textbf{F}_{i+1/2,j}(t)=\frac{1}{\Delta
			y}\int_{y_{j-1/2}}^{y_{j+1/2}}\textbf{F}_{i+1/2}(y,t)\text{d}y=\sum_{\ell=1}^2\omega_\ell \textbf{F}_{i+1/2,j_\ell}(t),
	\end{align}
	where $\omega_1=\omega_2=1/2$ are weights for the Gaussian
	quadrature point $\displaystyle y_{j_\ell}=y_j+\frac{(-1)^{\ell}}{2\sqrt{3}}\Delta y$, $\ell= 1, 2$,
	for a fourth-order accuracy.
	$\textbf{F}_{i+1/2,j_\ell}(t)$ are numerical fluxes and can be obtained
	as follows
	\begin{align}\label{flux-general}
		\textbf{F}_{i+1/2,j_\ell}(t)=\int \bm{\psi} u f(x_{i+1/2},y_{j_\ell},t,u,v,\xi)\text{d}\Xi,
	\end{align}
	where $f(x_{i+1/2},y_{j_\ell},t,u,v,\xi)$ is the gas distribution
	function at the cell interface. To construct the numerical
	fluxes, the integral solution of BGK equation Eq.\eqref{bgk} is used
	\begin{equation}\label{integral1}
		f(x_{i+1/2},y_{j_\ell},t,u,v,\xi)=\frac{1}{\tau}\int_0^t g(x',y',t',u,v,\xi)e^{-(t-t')/\tau}dt'\\
		+e^{-t/\tau}f_0(-ut,-vt,u,v,\xi),
	\end{equation}
	where $(x_{i+1/2}, y_{j_\ell})=(0,0)$ is the location for flux evaluation,
	and $x_{i+1/2}=x'+u(t-t')$ and $y_{j_\ell}=y'+v(t-t')$ are the trajectory of
	particles. $f_0$ is the initial gas distribution function, and $g$ is the corresponding
	equilibrium state. The integral solution states a physical process from the particle-free transport in $f_0$ in the kinetic scale
	to the hydrodynamic flow evolution in the integral of $g$ term.
	The flow behavior at the cell interface depends on the ratio of time step
	to the local particle collision time $\Delta t/\tau$.
	
	To construct a time evolution solution of a gas distribution function at a cell interface efficiently, a simplified third-order gas distribution proposed in \cite{zhou2017simplification} is adopted, which can be written as
	\begin{align}\label{3rd-simplify-flux}
		f(x_{i+1/2},y_{j_\ell},t,u,v,\xi)&=g^c+{A}^c g^ct+\frac{1}{2}{a}^c_{tt}g^ct^2\nonumber\\
		&-\tau[({a}^c_xu+{a}^c_yv+{A}^c)g_0+({a}^c_{xt}u+{a}^c_{yt}v+{a}^c_{tt})g_0t]\nonumber\\
		&-e^{-t/\tau}g^c[1-({a}^c_{x}u+{a}^c_{y}v)t]\nonumber\\
		&+e^{-t/\tau}g_l[1-(a_{xl}u+a_{yl}v)t]H(u)\nonumber\\
		&+e^{-t/\tau}g_r[1-(a_{xr}u+a_{yr}v)t](1-H(u)),
	\end{align}
	where the notations are introduced as follows
	\begin{align*}
		&a_x=(\partial g/\partial x)/g, a_y=(\partial g/\partial y)/g,
		a_t=A=(\partial g/\partial t)/g,\\
		&a_{xx}=g_{xx}/g,a_{xy}=g_{xy}/g,a_{yy}=g_{yy}/g,\\
		&a_{xt}=g_{xt}/g,a_{yt}=g_{yt}/g,a_{tt}=g_{tt}/g.
	\end{align*}
	The determination of these coefficients is simplified as
	\begin{align}\label{3rd-simplify-coe}
		\begin{cases}
			\displaystyle
			\left\langle a_{x}\right\rangle=\frac{\partial{\textbf{W}}}{\partial{x}},
			\left\langle a_{y}\right\rangle=\frac{\partial{\textbf{W}}}{\partial{y}},
			\left\langle a_{x}u+a_{y}v+a_{t}\right\rangle=0,
			\\
			\displaystyle
			\left\langle a_{xx}\right\rangle=\frac{\partial^2{\textbf{W}}}{\partial{x^2}},
			\left\langle a_{xy}\right\rangle=\frac{\partial^2{\textbf{W}}}{\partial{x}\partial{y}},
			\left\langle a_{yy}\right\rangle=\frac{\partial^2{\textbf{W}}}{\partial{y^2}},
			\\
			\displaystyle\left\langle a_{xx}u+a_{xy}v+a_{xt}\right\rangle=0,
			\\
			\displaystyle\left\langle a_{xy}u+a_{yy}v+a_{yt}\right\rangle=0,
			\\
			\displaystyle\left\langle a_{xt}u+a_{yt}v+a_{tt}\right\rangle=0,
		\end{cases}
	\end{align}
	where the superscripts and subscripts on these coefficients $a_{x},...,a_{tt}$ are
	omitted without ambiguity and $\left\langle ...\right\rangle$ are the moments of a gas distribution function defined by
	\begin{equation}
		\langle(\ldots)\rangle=\int \bm{\psi}(\ldots) g \mathrm{~d} \Xi.
	\end{equation}
	With the same third-order accuracy, the above
	simplified distribution function can speed up the flux calculation
	$4$ times in comparison to the complete gas distribution function in
	2-D case. The details of the conservative variables' reconstructions are shown in the next subsection.

	\subsection{Implementation of the spatial reconstruction}
	
	The  non-compact spatial reconstruction adopts the classical WENO-Z method and details are referred to \cite{ji2019high}.

	For compact spatial reconstruction, the Hermite WENO(HWENO) is introduced. HWENO was initially introduced to be one of the limiters of the RKDG method \cite{qiu2005hermite}. Varying from the classical WENO scheme, it utilized the derivative information of the numerical solutions and made it possible to achieve a compact reconstruction scheme, from which the compact spatial reconstruction of the gas-kinetic scheme was inspired.
	
\begin{figure}[!h]
		\centering
		\includegraphics[width=0.45\textwidth]{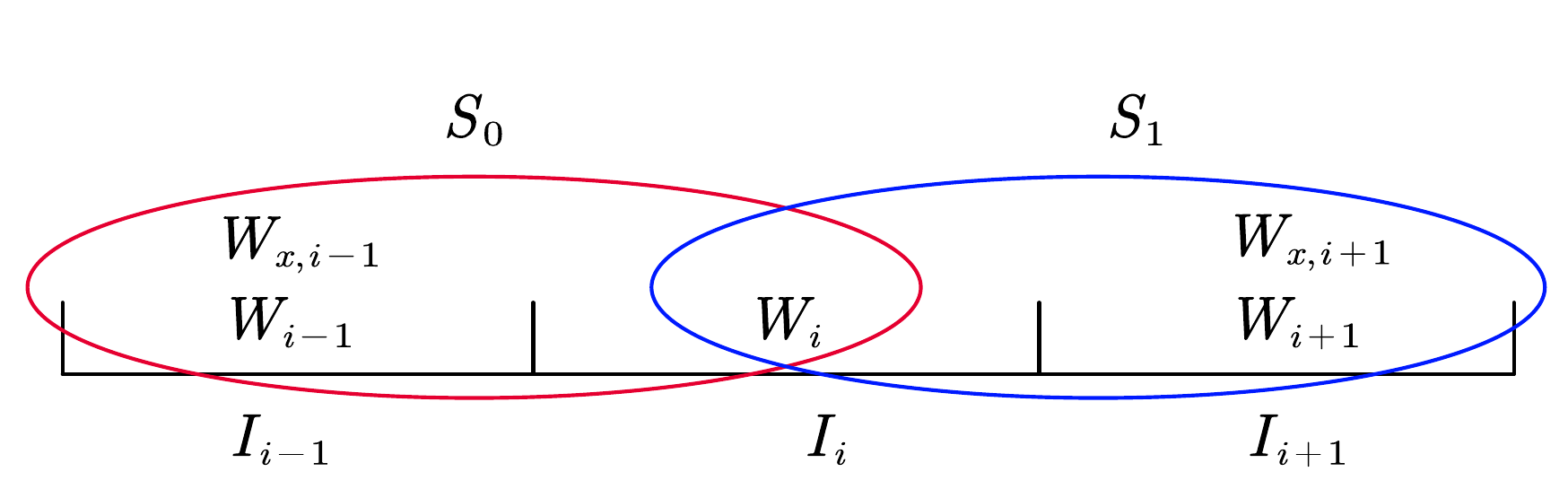}\hspace{5mm}
		\includegraphics[width=0.45\textwidth]{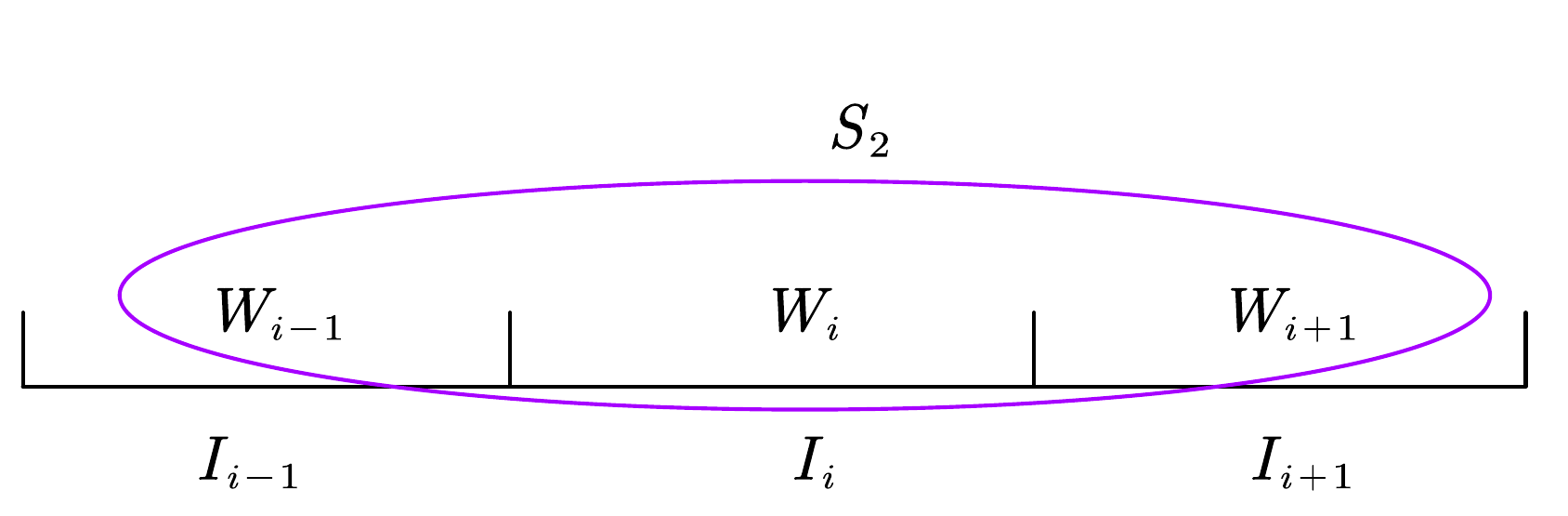}\hspace{5mm}
		\caption{\label{hweno-stencil-1d} Stencils for one dimensional compact gas-kinetic scheme reconstruction. }
	\end{figure}

	For one-dimensional cases, HWENO can be used to reconstruct the two sides of the interface value $W_{i+1/2}^l$ and $W_{i+1/2}^r$ at the interface $x_{i+1/2}$ \cite{qiu2004hermite}, where three sub-stencils are selected,
	\begin{equation}
		S_0=\left\{I_{i-1}, I_i\right\}, \quad S_1=\left\{I_i, I_{i+1}\right\}, \quad S_2=\left\{I_{i-1}, I_i, I_{i+1}\right\}.
	\end{equation}
	The quantities used for each stencil are shown in Figure \ref{hweno-stencil-1d}.
	Then, with the cell averaged value $W_i$ and the reconstructed values $W_{i\pm 1/2}$ on its two sides, a parabolic distribution of the conservative variables inside  $I_i$ is obtained to cater to the need for reconstruction of $f_0$. Next, the compatibility condition is used to obtain the initial equilibrium state $g_0$ with the reconstructed ${W}_{i+1/2}^l$ and ${W}_{i+1/2}^r$,
	\begin{equation}\label{compatibility2}
		\int \psi g^c \mathrm{~d} \Xi=W^c=\int_{u>0} \psi g_l \mathrm{~d} \Xi+\int_{u<0} \psi g_r \mathrm{~d} \Xi.
	\end{equation}
	
	To fully determine the slopes of the equilibrium state across the cell interface, the conservative variables across the cell interface are expanded as
	\begin{equation}
		{W}^c(x)=W_0+S_1\left(x-x_{i+1 / 2}\right)+\frac{1}{2} S_2\left(x-x_{i+1 / 2}\right)^2+\frac{1}{6} S_3\left(x-x_{i+1 / 2}\right)^3+\frac{1}{24} S_4\left(x-x_{i+1 / 2}\right)^4,
	\end{equation}
	where the derivatives are given by
\begin{equation}
	\left\{
	\begin{array}{l}
		{W}^c_x=S_1=\left[-\frac{1}{12}\left(W_{i+2}-W_{i-1}\right)+\frac{5}{4}\left(W_{i+1}-W_i\right)\right] / \Delta x \\
		\\
		{W}^c_{xx}=S_2=\left[-\frac{1}{8}\left(W_{i+2}+W_{i-1}\right)+\frac{31}{8}\left(W_{i+1}+W_i\right)-\frac{15}{2} W_0\right] / \Delta x^2.
	\end{array}
	\right.
\end{equation}

	Thus, the compact spatial reconstruction of the one-dimensional cases is fully given.

	For two-dimensional cases, the direction-by-direction reconstruction strategy \cite{zhang2011order} is adopted. For simplicity, the normal and tangential reconstruction procedure is illustrated briefly, more details are referred to \cite{ji2018compact}.
	
	For a clear description of the algorithm, we focus on the reconstruction based on the  $I_{i,j}$ whose x-direction index is $i$ and y-direction index is $j$.
	
	Step 1: Obtain the line-averaged non-central(left and right) values, $\widehat{W}_{i-1 / 2, j}^r$, $\widehat{W}_{i+1 / 2, j}^l$ as well as their normal derivatives $\widehat{W}_{x,i-1 / 2, j}^r$, $\widehat{W}_{x,i+1 / 2, j}^l$ and second-order derivatives $\widehat{W}_{xx,i-1 / 2, j}^r$, $\widehat{W}_{xx,i+1 / 2, j}^l$ by the cell-averaged values of conservative variables as well as their derivatives in a stencil centered on $I_{i,j}$, consisting of three adjacent cells $I_{i-1,j}$, $I_{i,j}$ and $I_{i+1,j}$, i.e., $\bar{W}_{i+k, j},k=0,\pm1$ and $\bar{W}_{x,i+k, j},k=\pm1$. In this step, a fifth-order WENO-Z-recipe HWENO reconstruction is adopted and the stencil is presented in Figure 2(a).
	
	Step 2: The central line averaged values $\widehat{W}_{i+1 / 2, j}^c$ are obtained by Eq.\eqref{compatibility2}. Then, the central line averaged normal derivatives  $\widehat{W}_{x,i+1 / 2, j}^c$ and second-order derivatives  $\widehat{W}_{xx,i+1 / 2, j}^c$ are obtained via fifth-order linear reconstruction with the five adjacent information, i.e.,  the cell averaged value $\bar{W}_{x,i+k, j},k=-1,0,1,2$ and the line-averaged values $\widehat{W}_{i+1 / 2, j}^c$. The reconstruction stencil is presented in Figure 2(b).
	
	Step 3: Obtain the values and their spatial derivatives at the Gaussian points, e.g. $\dot{W}_{i+1 / 2, j, 0}^{l,r}$, $\dot{W}_{y,i+1 / 2, j, k}^{l,r}$, $\dot{W}_{yy,i+1 / 2, j, k}^{l,r}$, $k=0,1$ with the reconstructed values $\widehat{W}_{i+1 / 2, j+k}^{l,r}$,$k=0,\pm1,\pm2$. Obtain the spatial derivatives at the Gaussian points, e.g. $\dot{W}_{x,i+1 / 2, j}$, $\dot{W}_{xy,i+1 / 2, j, k}^{l,r}$, $k=0,1$ with the reconstructed derivatives $\widehat{W}_{x,i+1 / 2, j+k}^{l,r},k=0,\pm1,\pm2$. Obtain the spatial derivatives at the Gaussian points, e.g. $\dot{W}_{xx,i+1 / 2, j, k}^{l,r},k=0,1$ with the reconstructed derivatives $\widehat{W}_{xx,i+1 / 2, j+k}^{l,r},k=0,\pm1,\pm2$. In this step, a fifth-order WENO reconstruction is adopted and the stencil is presented in Figure 2(c).
	
	Step 4: Obtain the central values and their spatial derivatives at the Gaussian points, i.e. $\dot{W}_{i+1 / 2, j, k}^c$, $\dot{W}_{y,i+1 / 2, j, k}^c$, $\dot{W}_{yy,i+1 / 2, j, k}^c,k=0,1$ with the values reconstructed $\widehat{W}_{i+1 / 2, j+k}^{c},k=0,\pm1,\pm2$. Obtain the central spatial derivatives at the Gaussian points, i.e. $\dot{W}_{x,i+1 / 2, j, k}^c$, $\dot{W}_{xy,i+1 / 2, j, k}^c,k=0,1$ with the derivatives reconstructed $\widehat{W}_{x,i+1 / 2, j+k}^{c},k=0,\pm1,\pm2$.
	Obtain the central spatial derivatives at the Gaussian points, i.e. $\dot{W}_{xx,i+1 / 2, j, k}^c,k=0,1$ with the derivatives reconstructed $\widehat{W}_{xx,i+1 / 2, j+k}^{c},k=0,\pm1,\pm2$.
	In this step, a fifth-order linear reconstruction is adopted and the stencil is the same as shown in Figure 2(c).
	
The reconstruction presented above is not entirely compact, as it incorporates a WENO5 reconstruction for the tangential part to ensure adequate accuracy. The complete stencil utilized for reconstructing values associated with inner sides of $I_{i,j}$ is depicted in Figure 2(d), with the dashed line portion specifically designated for tangential reconstruction.
	\color{black}{}
	
	\begin{figure}[!h]
		\centering
		\begin{subfigure}[b]{0.45\textwidth}
			\centering
			\includegraphics[width=\textwidth]{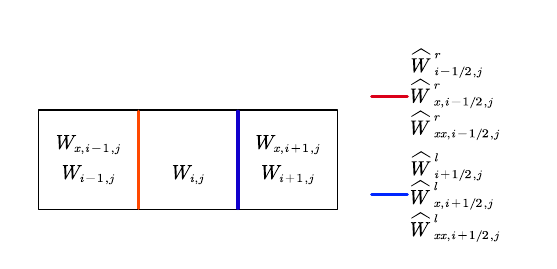}
			\caption{} % Caption for the first subfigure
		\end{subfigure}
		\hspace{5mm} % Horizontal space between subfigures
		\begin{subfigure}[b]{0.45\textwidth}
			\centering
			\includegraphics[width=\textwidth]{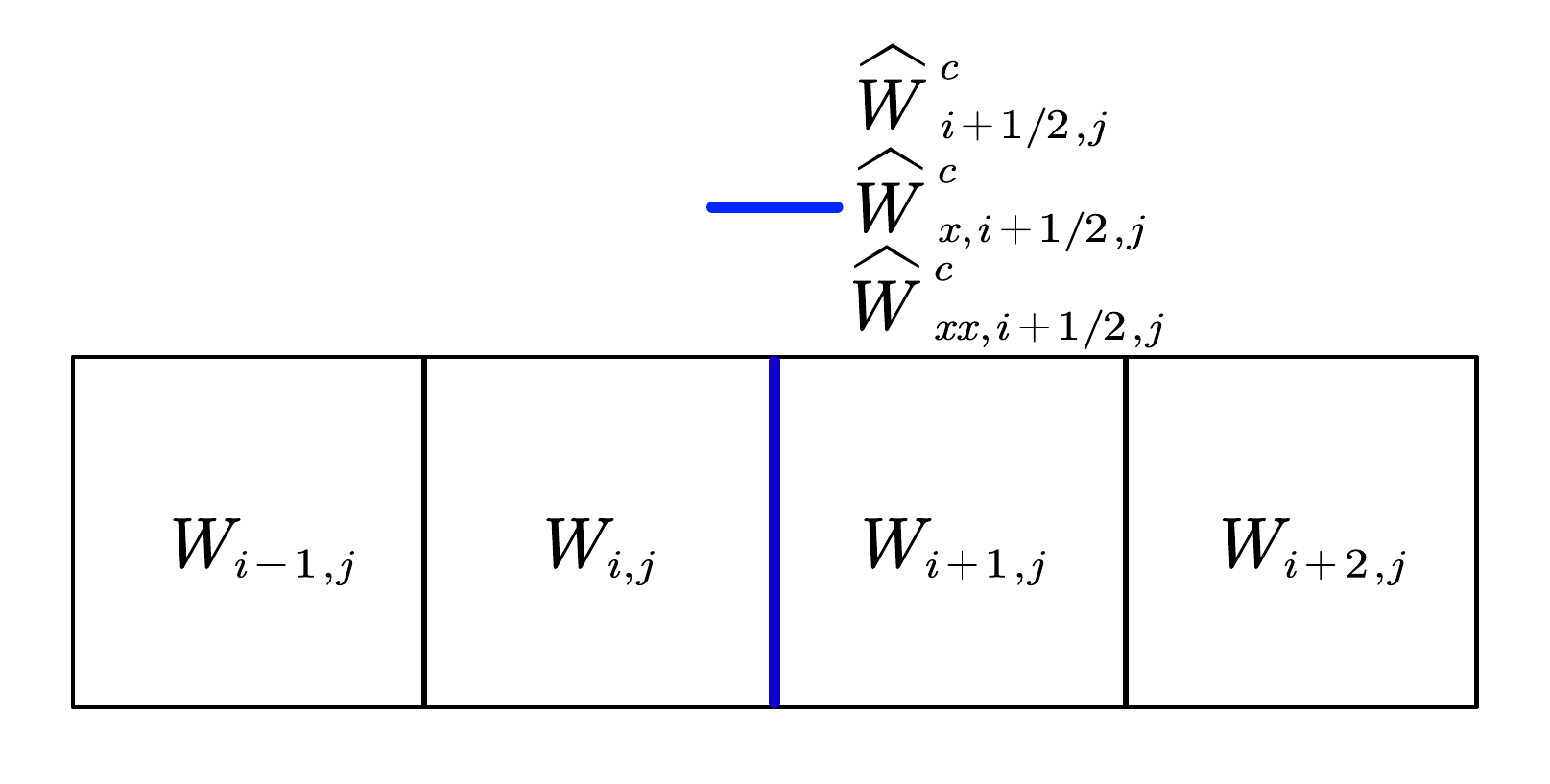}
			\caption{} % Caption for the second subfigure
		\end{subfigure}
		\begin{subfigure}[b]{0.45\textwidth}
			\centering
			\includegraphics[height=\textwidth]{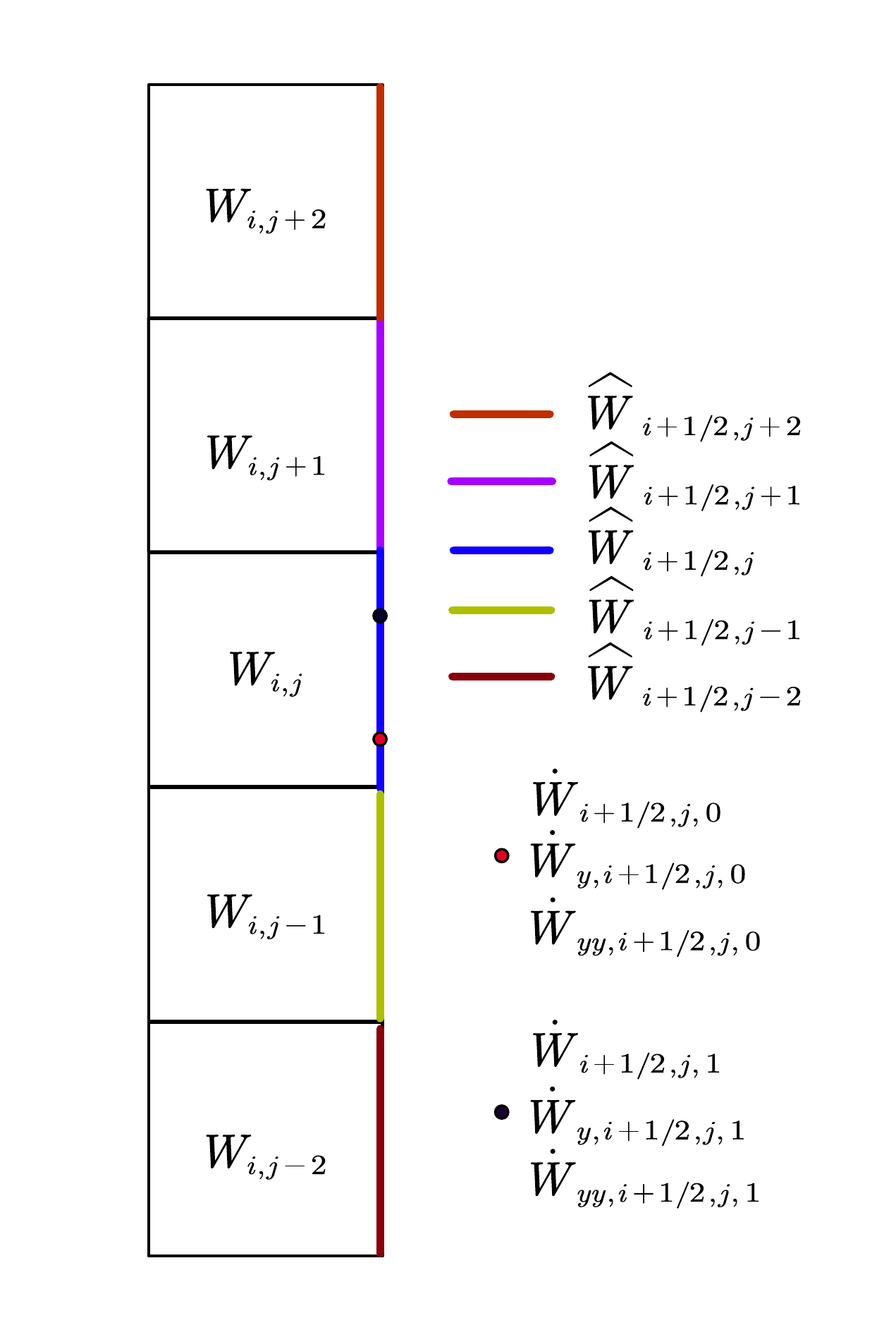}
			\caption{} % Caption for the second subfigure
		\end{subfigure}
		\hspace{5mm}
		\begin{subfigure}[b]{0.45\textwidth}
			\centering
			\includegraphics[height=\textwidth]{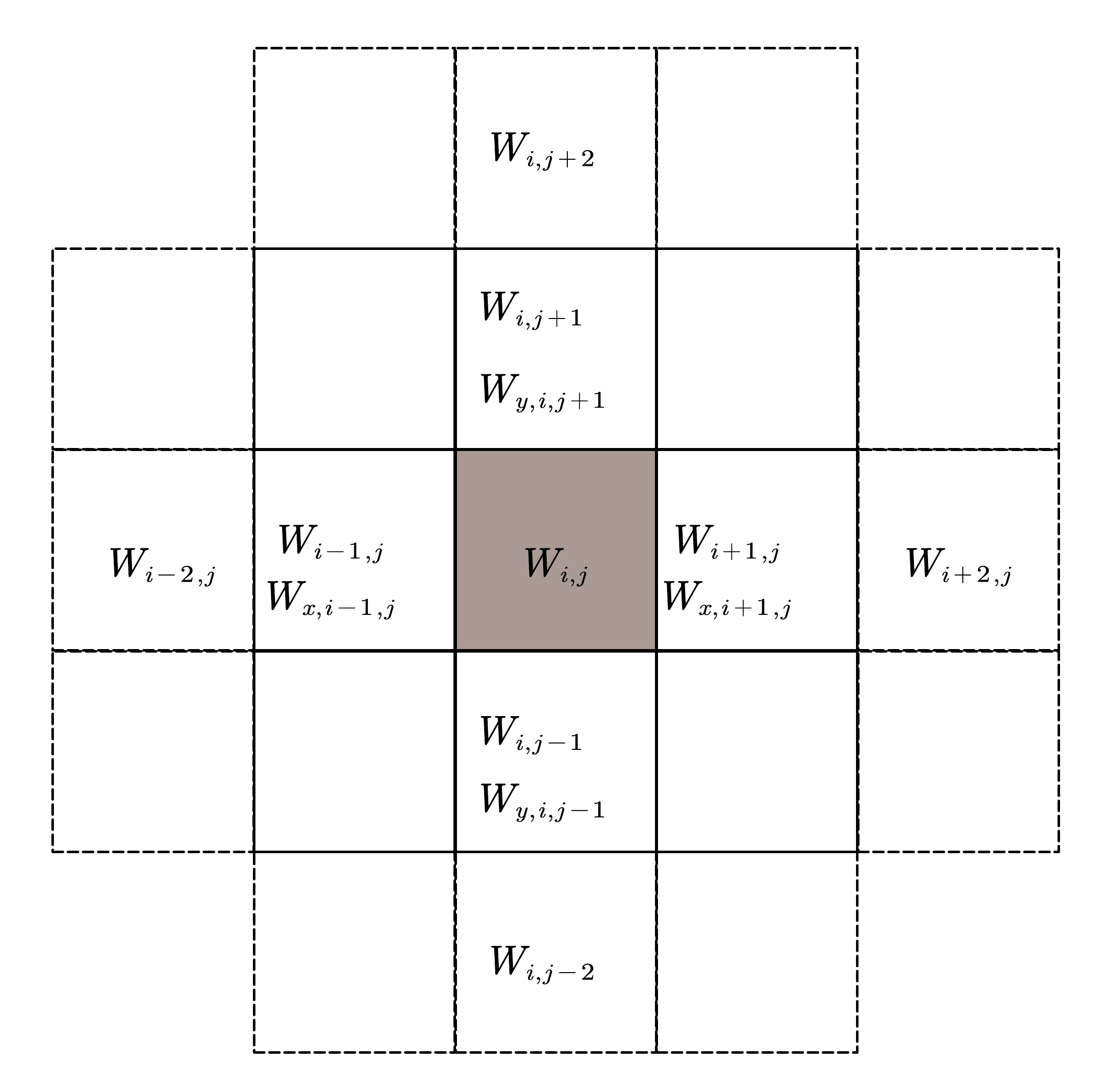}
			\caption{} % Caption for the second subfigure
		\end{subfigure}
		
		\caption{Stencils for two-dimensional compact gas-kinetic scheme reconstruction.}
		\label{hweno-stencil}
	\end{figure}

	To resolve the problems with shock waves, characteristic variables are used, consistent with the RKDG method. In CGKS, three vectors are stored in one specific cell in each time step, $W_{ij}$,$W_{x,ij}$ and $W_{y,ij}$, namely cell-averaged values and two-dimensional cell-averaged spatial derivatives respectively.

	\subsection{Two-stage fourth-order temporal discretization}
	The two-stage fourth-order temporal discretization \cite{li2016two} is based on the semi-discrete finite volume scheme, which is written as
	\begin{align*}
		\frac{\text{d} {W}_{ij}}{\text{d} t}=-\frac{1}{\Delta x}
		(F_{i+1/2,j}(t)-F_{i-1/2,j}(t))-\frac{1}{\Delta y}
		(G_{i,j+1/2}(t)-G_{i,j-1/2}(t)):=\mathcal{L}(W_{ij}),
	\end{align*}
	where $\mathcal{L}$ is the numerical operator for spatial derivative of flux, $F$ and $G$ are obtained by the Gaussian quadrature in Eq.\eqref{gauss-flux}.
	By adopting techniques in \cite{pan2016efficient}, this kind of temporal discretization is applied to the gas-kinetic scheme.
	
	The foundation of this work is to introduce the time derivatives of the operator for the spatial derivative and take advantage of the temporal accuracy of the gas-kinetic flux.
	
	Consider the temporal evolution at $t_n$, the Taylor expansion in time can be written as
	\begin{equation}
		W^{n+1}=W^n+\Delta t\mathcal
		{L}(W^n)+\frac{1}{6}\Delta t^2\frac{\partial}{\partial
			t}\mathcal{L}(W^n)+\frac{1}{3}\Delta t^2\frac{\partial}{\partial
			t}\mathcal{L}(W^*),
	\end{equation}
	to achieve fourth-order accuracy in temporal evolution, the intermediate state is introduced at half of the time step \cite{li2016two}, which holds the expansion simultaneously:
	\begin{equation}
		W^{*}=W^n+\frac{1}{2}\Delta t\mathcal
		{L}(W^n)+\frac{1}{8}\Delta t^2\frac{\partial}{\partial
			t}\mathcal{L}(W^n),
	\end{equation}

	To obtain the complete temporal scheme, the temporal derivatives of fluxes ought to be calculated, which are easy to implement, since the gas-kinetic flux has temporal third-order accuracy.
	In practice, the temporal derivatives of fluxes can be derived by the linear approximation,
	\begin{equation}\label{second}
		\begin{aligned}
			F_{i+1/2,j}(W^n,t_n)&=(4\mathbb{F}_{i+1/2,j}(W^n,\Delta t/2)-\mathbb{F}_{i+1/2,j}(W^n,\Delta t))/\Delta t,\\
			\partial_t F_{i+1/2,j}(W^n,t_n)&=4(\mathbb{F}_{i+1/2,j}(W^n,\Delta t)-2\mathbb{F}_{i+1/2,j}(W^n,\Delta t/2))/\Delta
			t^2,
		\end{aligned}
	\end{equation}
	where the notation of the Integral $\mathbb{F}_{i+1/2,j}(W^n,\delta)$ is defined as
	\begin{align*}
		\mathbb{F}_{i+1/2,j}(W^n,\delta)
		=\int_{t_n}^{t_n+\delta}F_{i+1/2,j}(W^n,t)dt&=\sum_{\ell=1}^2\omega_\ell \int_{t_n}^{t_n+\delta}\int
		u \psi f(x_{i+1/2},y_{j_\ell},t,u, v,\xi)\text{d}u\text{d}v\text{d}\xi\text{d}t.
	\end{align*}
	
	Similarly, the numerical fluxes $G_{i,j+1/2}(W^n,t)$ in the y-direction can also be derived from the algorithm above. Thus, the temporal derivatives of the spacial operator can be written as
	\begin{align}
		\frac{\partial}{\partial t}\mathcal{L}(W_{ij}^n)&= - \frac{1}{\Delta x}(\partial_t F_{i+1/2,j}(W^n,t_n)-\partial_t F_{i-1/2,j}(W^n,t_n)) \nonumber \\
		& - \frac{1}{\Delta y}(\partial_t G_{i,j+1/2}(W^n,t_n)-\partial_t G_{i,j-1/2}(W^n,t_n))\label{operator-2}.
	\end{align}
	
	The point-wise values at a cell interface can be obtained by taking moments of the time-dependent distribution function
	in Eq.\eqref{3rd-simplify-flux}
	\color{black}{}
	\begin{align}\label{point}
		W_{i+1/2,j_\ell}(t)=\int\psi f(x_{i+1/2},y_{j_\ell},t,u,v,\xi)\text{d}\Xi.
	\end{align}
	
	The evolution of the gas distribution function with similar techniques \cite{ji2018compact} is obtained.

	To utilize the two-stage, fourth-order temporal
	discretization for the gas distribution function, the gas distribution function should be approximated by a quadratic polynomial,
	\begin{align*}
		f(t)=f(x_{i+1/2},y_{j_\ell},t,u, v,\xi)=f^n+
		f_{t}^n(t-t^n)+\frac{1}{2}f_{tt}^n(t-t^n)^2.
	\end{align*}
	the coefficients $f^n, f_{t}^n$ and $f_{tt}^n$ can be
	determined
	\begin{align*}
		f^n&=f(0),\\
		f^n_t&=(4f(\Delta t/2)-3f(0)-f(\Delta t))/\Delta t,\\
		f^n_{tt}&=4(f(\Delta t)+f(0)-2f(\Delta t/2))/\Delta t^2.
	\end{align*}
	Thus, $f^*$ and $f^{n+1}$ are entirely determined at the cell interface to evaluate macroscopic flow variables.

	\section{Numerical examples}
	
	For the RKDG method, the CFL number of all cases is taken as 0.18 for third-order($P^{2})$, and 0.08 for fourth-order($P^{4})$ both in one-dimensional and two-dimensional cases. Third-order Runge-Kutta method is adopted for  evolution if not mentioned specifically.  As shown in \cite{qiu2006numerical}, the Lax–Friedrichs(LF) flux costs the least CPU time, and the Harten–Lax–van Leer(HLL) flux generally performs well, taking both accuracy and efficiency into account. In this section, numerical tests will be presented to make comparisons between the RKDG method with the gas-kinetic scheme. All simulations are for compressible inviscid flow.
	Two numerical fluxes are used for the RKDG method in our work. However, there is little difference in the performance of various numerical fluxes on two-dimensional discontinuity problems. Thus, only LF flux is used in the double Mach reflection problem.
	
	For the gas-kinetic scheme, the artificial collision time $\tau$ for inviscid flow is taken as
	\begin{align*}
		\tau=\epsilon \Delta t+C\displaystyle|\frac{p_l-p_r}{p_l+p_r}|\Delta
		t,
	\end{align*}
	where $\epsilon=0.01$ and $C=1$.
	
	In our work, both the non-compact gas-kinetic scheme(GKS) in third-order and fifth-order, as well as the fifth-order compact gas-kinetic scheme(CGKS), are tested in order to make the comparison more comprehensive.
	At the same time, in order to ensure consistency with the multi-resolution limiters, the WENO-Z recipe is used to measure the smoothness indicator factor and nonlinear weights. The non-compact gas-kinetic scheme(GKS) for the third-order is based on WENO3-type spatial reconstruction and for the fifth-order, WENO5-type spatial reconstruction is used. The compact gas-kinetic scheme (CGKS), adopts the fifth-order HWENO-Z reconstruction and a fourth-order in temporal reconstruction. Both the non-compact and compact GKS use two-stage fourth-order time-stepping method.
	
	The efficiency comparison cases are performed on CPU 12th Gen Intel(R) Core(TM) i7-12700H @2.30GHz with one single processor running. The theoretical bandwidth of memory is 38.4GB/s. data structures of both two methods are based on Array of Structures(AOS). All the numerical results are obtained by our in-house C++ solver.

	\subsection{One dimensional accuracy tests}
	
	In this subsection, one dimensional advection of density perturbation is tested, and the initial
	condition is given as follows
	\begin{align*}
		\rho(x)=1+0.2\sin(\pi x),\ \  U(x)=1,\ \ \  p(x)=1, x\in[-1,1].
	\end{align*}
	With the periodic boundary condition, and the analytic
	solution is
	\begin{align*}
		\rho(x,t)=1+0.2\sin(\pi(x-t)),\ \ \  U(x,t)=1,\ \ \  p(x,t)=1.
	\end{align*}
	In this computation, the computational domain is partitioned by uniform meshes, and the final time is $t=2s$. And the classical fourth-order Runge-Kutta method  \cite{butcher1996history} is adopted for RKDG-$P^{4}$ method to guarantee that the spatial error dominates.
	The $L^{1}$ and $L^{2}$ errors, along with their orders, are presented in Table 1-2 for the third-order RKDG and GKS, and in Table 3-4 for the fifth-order. And we measure the efficiency of the method using the CPU time calculated from the $L^{1}$ and $L^{2}$ errors, which is shown in Figure 2 and Figure 3.

	\begin{table}[]
		\centering
		\begin{tabular}{|l|llll|llll|}
			\hline
			$N_{cell}$ & \multicolumn{4}{c|}{RKDG-$P^{2}$ with LF flux}                                                               & \multicolumn{4}{c|}{RKDG-$P^{2}$ with HLL flux}                                                              \\ \hline
			& \multicolumn{1}{l|}{$L^{1}$ error} & \multicolumn{1}{l|}{order}  & \multicolumn{1}{l|}{$L^{2}$  error} & order  & \multicolumn{1}{l|}{$L^{1}$ error} & \multicolumn{1}{l|}{order}  & \multicolumn{1}{l|}{$L^{2}$  error} & order  \\ \hline
			20   & \multicolumn{1}{l|}{3.00E-05} & \multicolumn{1}{l|}{}       & \multicolumn{1}{l|}{3.35E-05} &        & \multicolumn{1}{l|}{1.72E-05} & \multicolumn{1}{l|}{}       & \multicolumn{1}{l|}{1.93E-05} &        \\ \hline
			40   & \multicolumn{1}{l|}{3.93E-06} & \multicolumn{1}{l|}{2.93} & \multicolumn{1}{l|}{4.39E-06} & 2.93 & \multicolumn{1}{l|}{2.17E-06} & \multicolumn{1}{l|}{2.99} & \multicolumn{1}{l|}{2.43E-06} & 3.00 \\ \hline
			80   & \multicolumn{1}{l|}{4.98E-07} & \multicolumn{1}{l|}{2.98} & \multicolumn{1}{l|}{5.56E-07} & 3.00 & \multicolumn{1}{l|}{2.72E-07} & \multicolumn{1}{l|}{3.00} & \multicolumn{1}{l|}{3.04E-07} & 3.00 \\ \hline
			160  & \multicolumn{1}{l|}{6.25E-08} & \multicolumn{1}{l|}{3.00} & \multicolumn{1}{l|}{6.98E-08} & 3.00 & \multicolumn{1}{l|}{3.40E-08} & \multicolumn{1}{l|}{3.00} & \multicolumn{1}{l|}{3.80E-08} & 3.00 \\ \hline
			320  & \multicolumn{1}{l|}{7.81E-09} & \multicolumn{1}{l|}{3.00} & \multicolumn{1}{l|}{8.73E-09} & 3.00 & \multicolumn{1}{l|}{4.25E-09} & \multicolumn{1}{l|}{3.00} & \multicolumn{1}{l|}{4.75E-09} & 3.00 \\ \hline
			640  & \multicolumn{1}{l|}{9.77E-10} & \multicolumn{1}{l|}{3.00} & \multicolumn{1}{l|}{1.09E-09} & 3.00 & \multicolumn{1}{l|}{5.31E-10} & \multicolumn{1}{l|}{3.00} & \multicolumn{1}{l|}{5.93E-10} & 3.00 \\ \hline
		\end{tabular}
		\caption{One dimensional advection of density perturbation: accuracy test for classical RKDG-$P^{2}$ method without limiters.}
	\end{table}

	% Please add the following required packages to your document preamble:
	% \usepackage{graphicx}
	\begin{table}[]
		\centering
		\begin{tabular}{|l|llll|}
			\hline
			$N_{cell}$ & \multicolumn{4}{c|}{Third-order GKS}                                                                             \\ \hline
			& \multicolumn{1}{l|}{$L^{1}$ error} & \multicolumn{1}{l|}{order}  & \multicolumn{1}{l|}{$L^{2}$ error} & order  \\ \hline
			20   & \multicolumn{1}{l|}{9.29E-04} & \multicolumn{1}{l|}{}       & \multicolumn{1}{l|}{1.03E-03} &        \\ \hline
			40   & \multicolumn{1}{l|}{1.17E-04} & \multicolumn{1}{l|}{2.99} & \multicolumn{1}{l|}{1.30E-04} & 2.99 \\ \hline
			80   & \multicolumn{1}{l|}{1.47E-05} & \multicolumn{1}{l|}{3.00} & \multicolumn{1}{l|}{1.62E-05} & 3.00 \\ \hline
			160  & \multicolumn{1}{l|}{1.83E-06} & \multicolumn{1}{l|}{3.00} & \multicolumn{1}{l|}{2.03E-06} & 3.00 \\ \hline
			320  & \multicolumn{1}{l|}{2.29E-07} & \multicolumn{1}{l|}{3.00} & \multicolumn{1}{l|}{2.53E-07} & 3.00 \\ \hline
			640  & \multicolumn{1}{l|}{2.86E-08} & \multicolumn{1}{l|}{3.00} & \multicolumn{1}{l|}{3.16E-08} & 3.00 \\ \hline
		\end{tabular}
		\caption{One dimensional advection of density perturbation: accuracy test for third-order GKS at smooth reconstruction.}
	\end{table}

	\begin{table}[]
		\centering
		\begin{tabular}{|l|llll|llll|}
			\hline
			$N_{cell}$ & \multicolumn{4}{c|}{RKDG-$P^{4}$ with LF flux}                                                               & \multicolumn{4}{c|}{RKDG-$P^{4}$ with HLL flux  }                                                            \\ \hline
			& \multicolumn{1}{l|}{$L^{1}$ error} & \multicolumn{1}{l|}{order}  & \multicolumn{1}{l|}{$L^{2}$ error} & order  & \multicolumn{1}{l|}{$L^{1}$ error} & \multicolumn{1}{l|}{order}  & \multicolumn{1}{l|}{$L^{2}$ error} & order  \\ \hline
			10   & \multicolumn{1}{l|}{2.39E-07} & \multicolumn{1}{l|}{}       & \multicolumn{1}{l|}{2.69E-07} &        & \multicolumn{1}{l|}{1.56E-07} & \multicolumn{1}{l|}{}       & \multicolumn{1}{l|}{1.73E-07} &        \\ \hline
			20   & \multicolumn{1}{l|}{8.76E-09} & \multicolumn{1}{l|}{4.77} & \multicolumn{1}{l|}{9.79E-09} & 4.78 & \multicolumn{1}{l|}{4.98E-09} & \multicolumn{1}{l|}{4.97} & \multicolumn{1}{l|}{5.57E-09} & 4.95 \\ \hline
			40   & \multicolumn{1}{l|}{2.87E-10} & \multicolumn{1}{l|}{4.93} & \multicolumn{1}{l|}{3.21E-10} & 4.93 & \multicolumn{1}{l|}{1.58E-10} & \multicolumn{1}{l|}{4.98} & \multicolumn{1}{l|}{1.77E-10} & 4.98 \\ \hline
			80   & \multicolumn{1}{l|}{9.06E-12} & \multicolumn{1}{l|}{4.98} & \multicolumn{1}{l|}{1.02E-11} & 4.98 & \multicolumn{1}{l|}{4.95E-12} & \multicolumn{1}{l|}{4.99} & \multicolumn{1}{l|}{5.55E-12} & 4.99 \\ \hline
			160  & \multicolumn{1}{l|}{2.85E-13} & \multicolumn{1}{l|}{4.99} & \multicolumn{1}{l|}{3.19E-13} & 4.99 & \multicolumn{1}{l|}{1.55E-13} & \multicolumn{1}{l|}{4.99} & \multicolumn{1}{l|}{1.75E-13} & 4.99 \\ \hline
		\end{tabular}
		\caption{One dimensional advection of density perturbation: accuracy test for classical RKDG-$P^{4}$ method without limiters.}
	\end{table}

	\begin{table}[]
		\centering
		\begin{tabular}{|l|llll|llll|}
			\hline
			$N_{cell}$ & \multicolumn{4}{c|}{Fifth-order GKS}                                                                   & \multicolumn{4}{c|}{Fifth-order CGKS}                                                                \\ \hline
			& \multicolumn{1}{l|}{L1 error}   & \multicolumn{1}{l|}{order}  & \multicolumn{1}{l|}{L2 error} & order  & \multicolumn{1}{l|}{L1 error} & \multicolumn{1}{l|}{order}  & \multicolumn{1}{l|}{L2 error} & order  \\ \hline
			10   & \multicolumn{1}{l|}{8.49E-04}   & \multicolumn{1}{l|}{}       & \multicolumn{1}{l|}{9.58E-04} &        & \multicolumn{1}{l|}{1.53E-04} & \multicolumn{1}{l|}{}       & \multicolumn{1}{l|}{1.72E-04} &        \\ \hline
			20   & \multicolumn{1}{l|}{2.81E-05}   & \multicolumn{1}{l|}{4.92} & \multicolumn{1}{l|}{3.12E-05} & 4.94 & \multicolumn{1}{l|}{4.80E-06} & \multicolumn{1}{l|}{4.99} & \multicolumn{1}{l|}{5.33E-06} & 5.01 \\ \hline
			40   & \multicolumn{1}{l|}{8.89E-07}   & \multicolumn{1}{l|}{4.98} & \multicolumn{1}{l|}{9.84E-07} & 4.99 & \multicolumn{1}{l|}{1.50E-07} & \multicolumn{1}{l|}{5.00} & \multicolumn{1}{l|}{1.66E-07} & 5.00 \\ \hline
			80   & \multicolumn{1}{l|}{2.78E-08}   & \multicolumn{1}{l|}{5.00} & \multicolumn{1}{l|}{3.08E-08} & 5.00 & \multicolumn{1}{l|}{4.73E-09} & \multicolumn{1}{l|}{4.98} & \multicolumn{1}{l|}{5.25E-09} & 4.98 \\ \hline
			160  & \multicolumn{1}{l|}{8.7091E-10} & \multicolumn{1}{l|}{5.00} & \multicolumn{1}{l|}{9.65E-10} & 5.00 & \multicolumn{1}{l|}{1.54E-10} & \multicolumn{1}{l|}{4.94} & \multicolumn{1}{l|}{1.72E-10} & 4.93 \\ \hline
		\end{tabular}
		\caption{One dimensional advection of density perturbation: accuracy test for fifth-order GKS and CGKS  at smooth reconstruction.}
		
	\end{table}

	\begin{figure}[!h]
		\centering
		\includegraphics[width=0.49\textwidth]{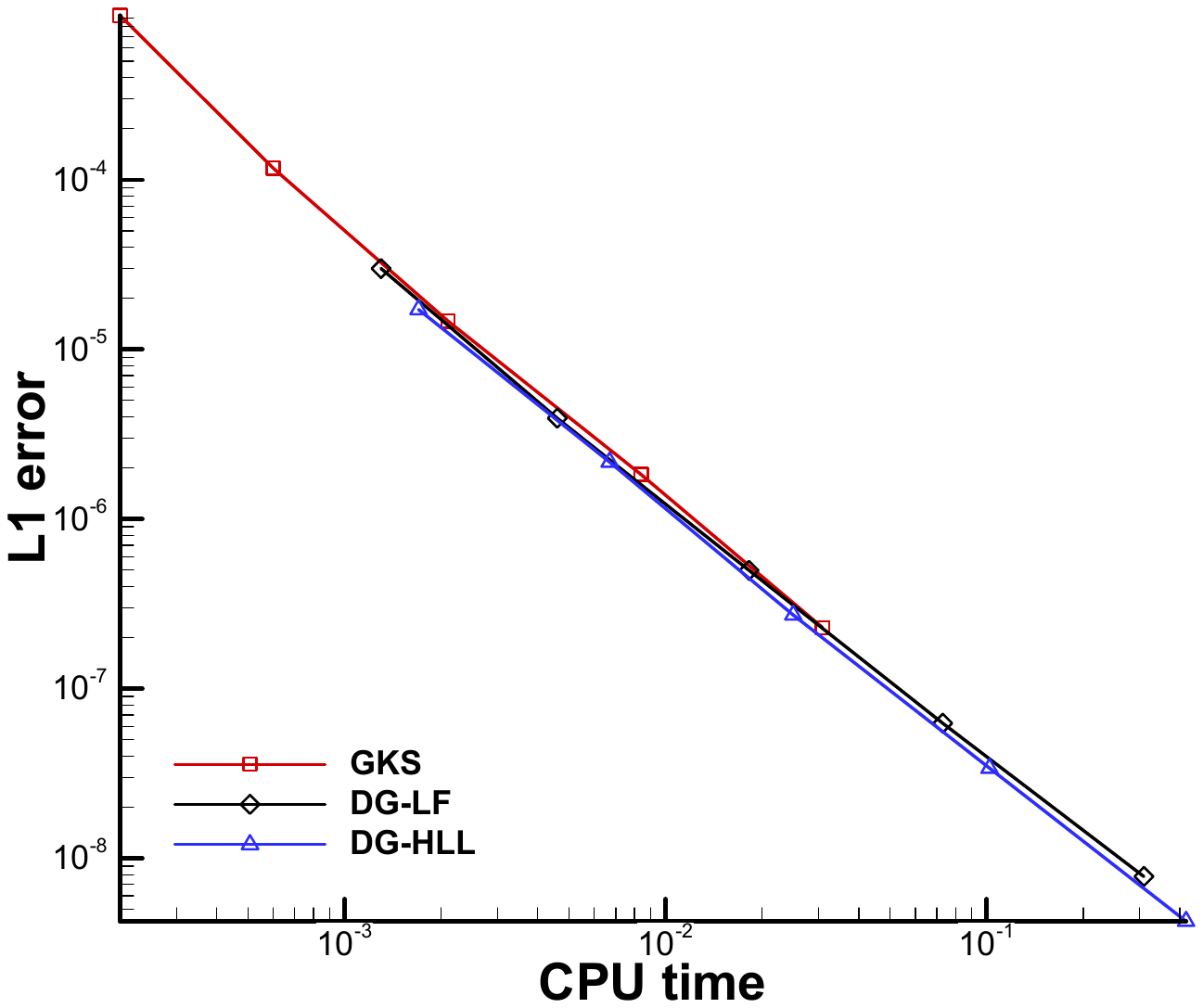}
		\includegraphics[width=0.49\textwidth]{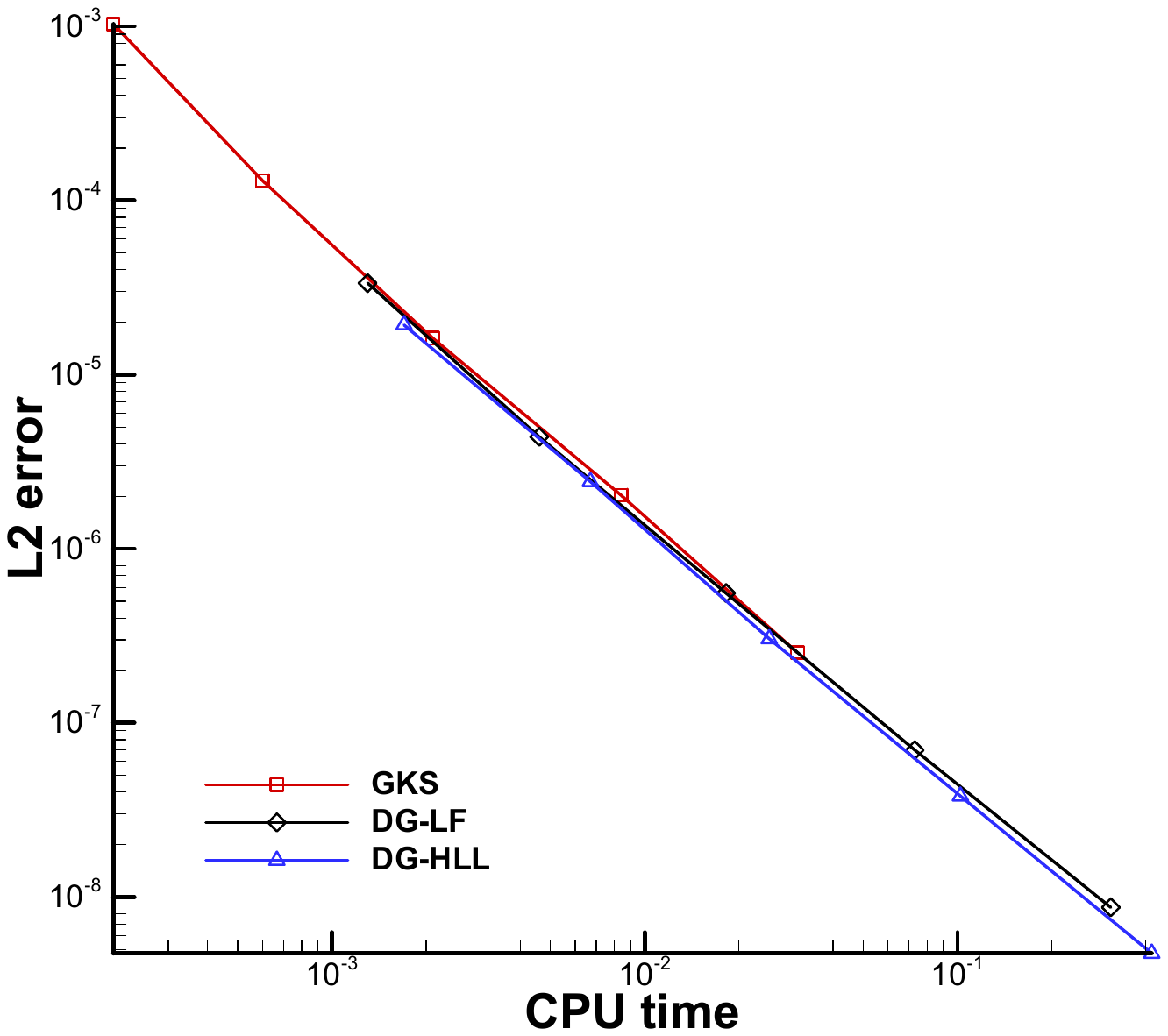}
		\caption{\label{accuracy-3-rd} Efficiency of simulation for smooth regions: comparisons between RKDG-$P^2$ method and third-order GKS method. The left one is based on $L^{1}$ error and the right one is based on $L^{2}$ error.}
	\end{figure}

	\begin{figure}[!h]
		\centering
		\includegraphics[width=0.49\textwidth]{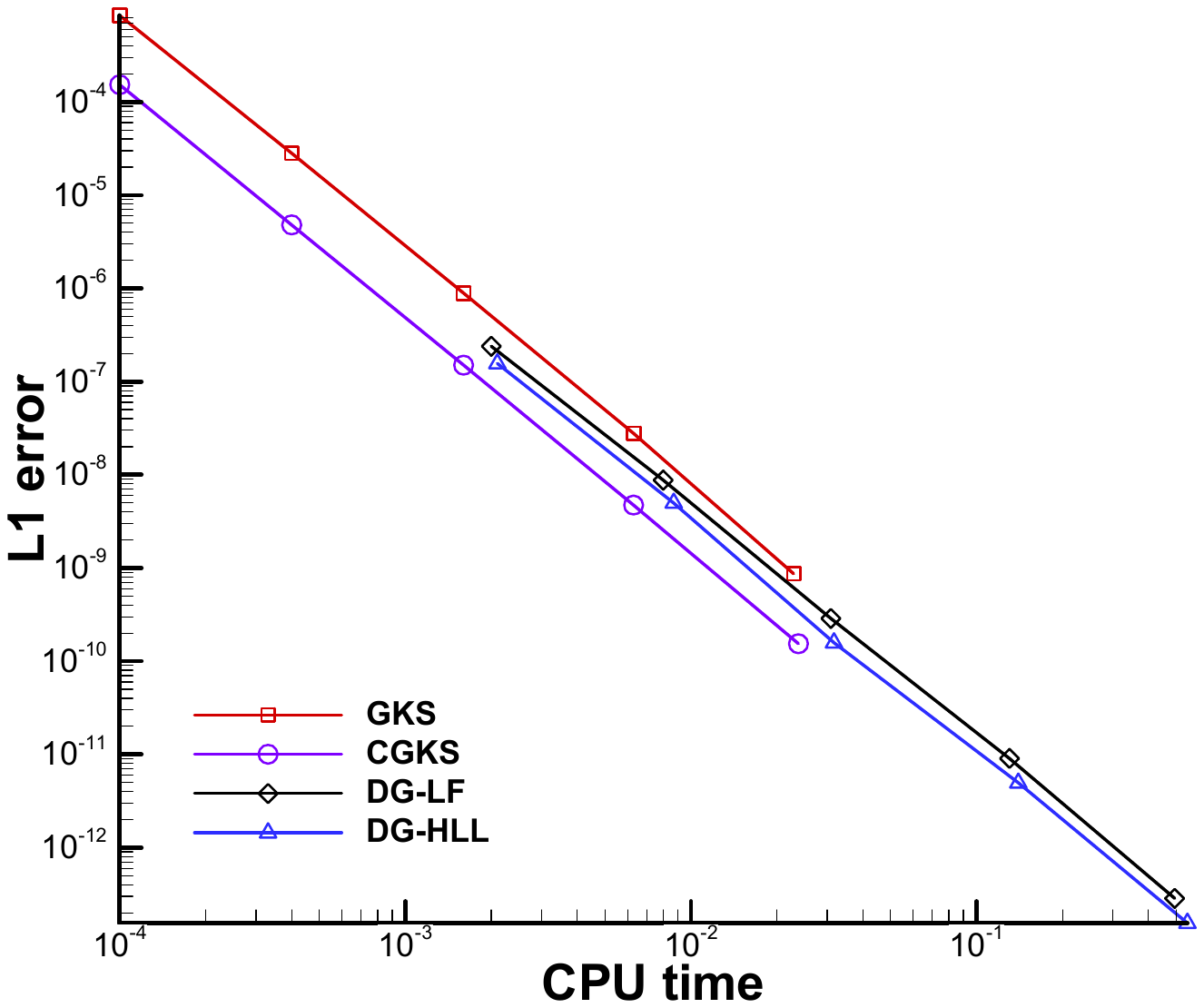}
		\includegraphics[width=0.49\textwidth]{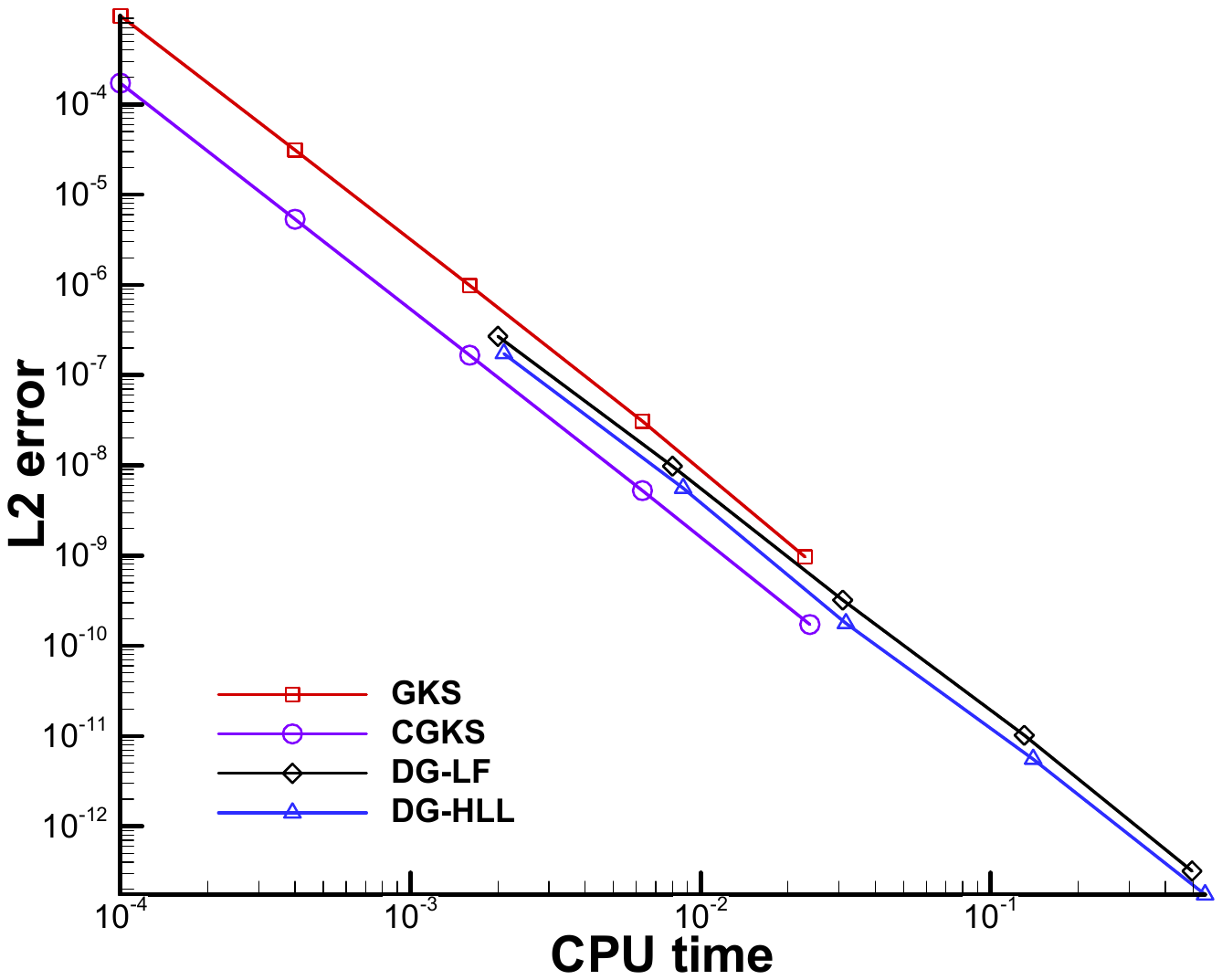}
		\caption{\label{accuracy-3-rg} Efficiency of simulation for smooth regions: comparisons between RKDG-$P^4$  and fifth-order GKS and CGKS. The left one is based on $L^{1}$ error and the right one is based on $L^{2}$ error.}
	\end{figure}
	
	According to the comparison of efficiency diagrams, it is apparent to see that in both third-order and fifth-order cases, the RKDG method with either LF or HLL flux has less error than the non-compact GKS  in the smooth region. However, the  CGKS  is shown to perform better than the RKDG-$P^{4}$ method with LF flux and HLL flux.
	
	\begin{figure}[!h]
		\centering
		\includegraphics[width=0.32\textwidth, trim=20mm 0mm 30mm 10mm, clip]{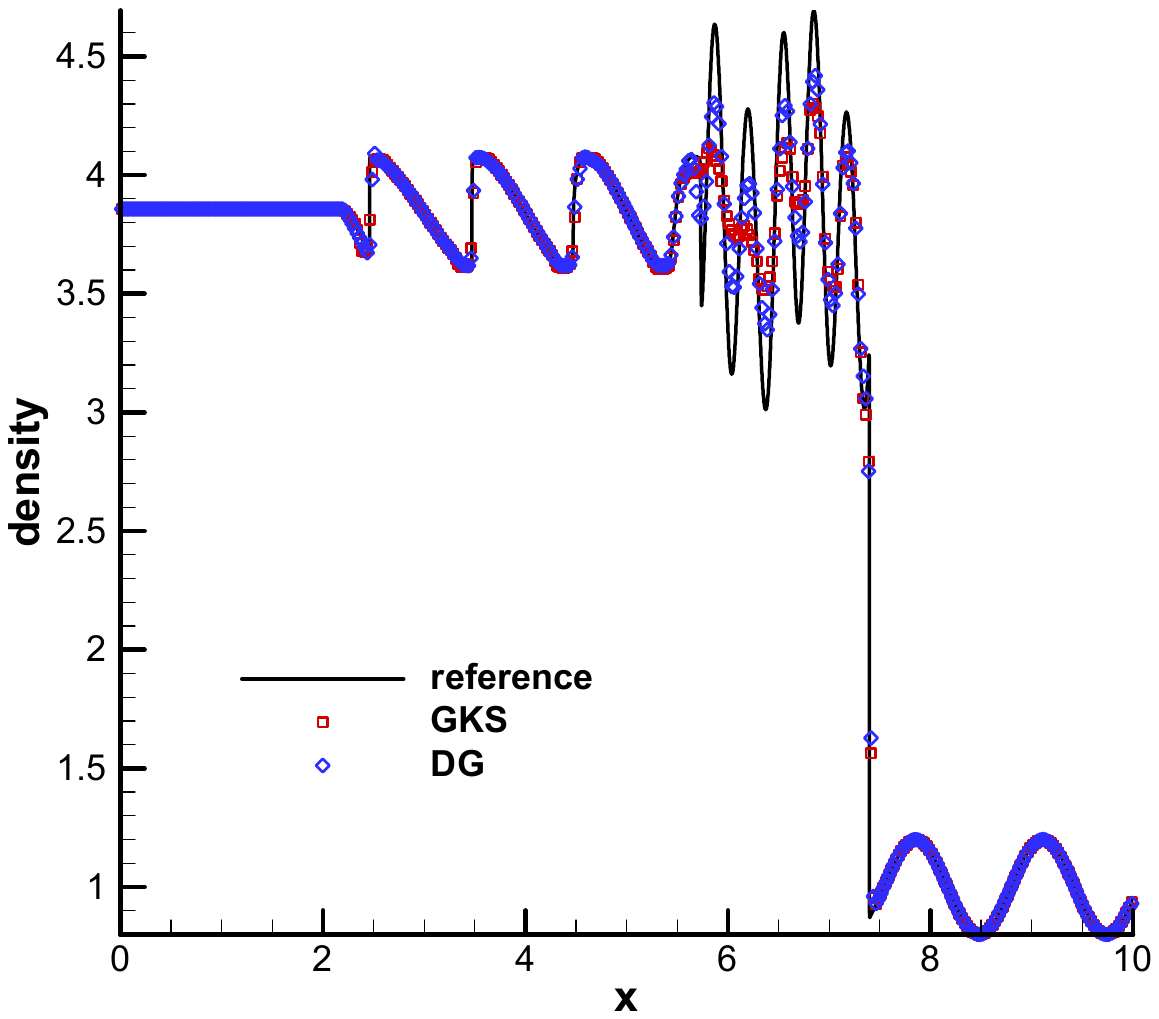}
		\includegraphics[width=0.32\textwidth, trim=20mm 0mm 30mm 10mm, clip]{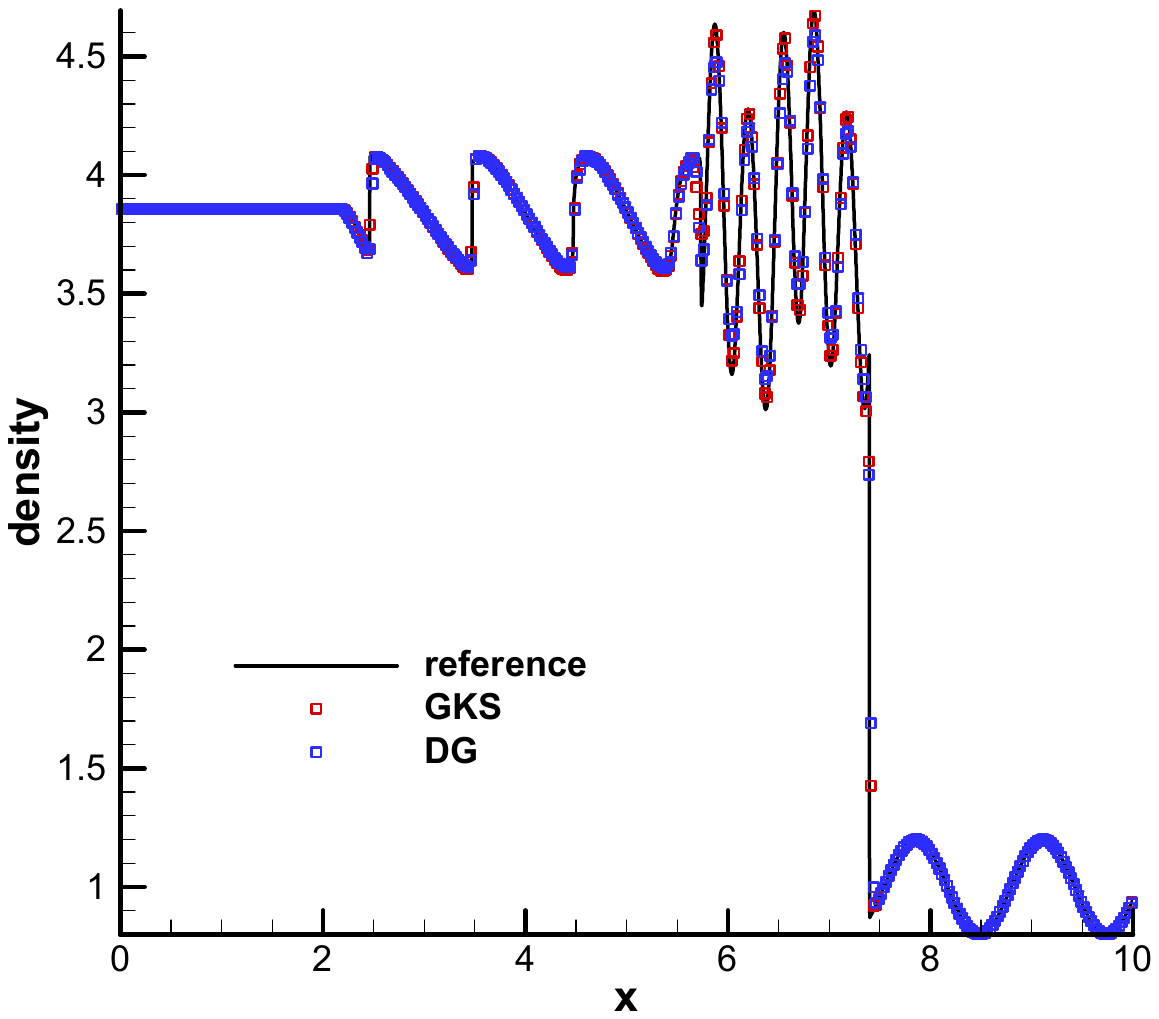}
		\includegraphics[width=0.32\textwidth, trim=20mm 0mm 30mm 10mm, clip]{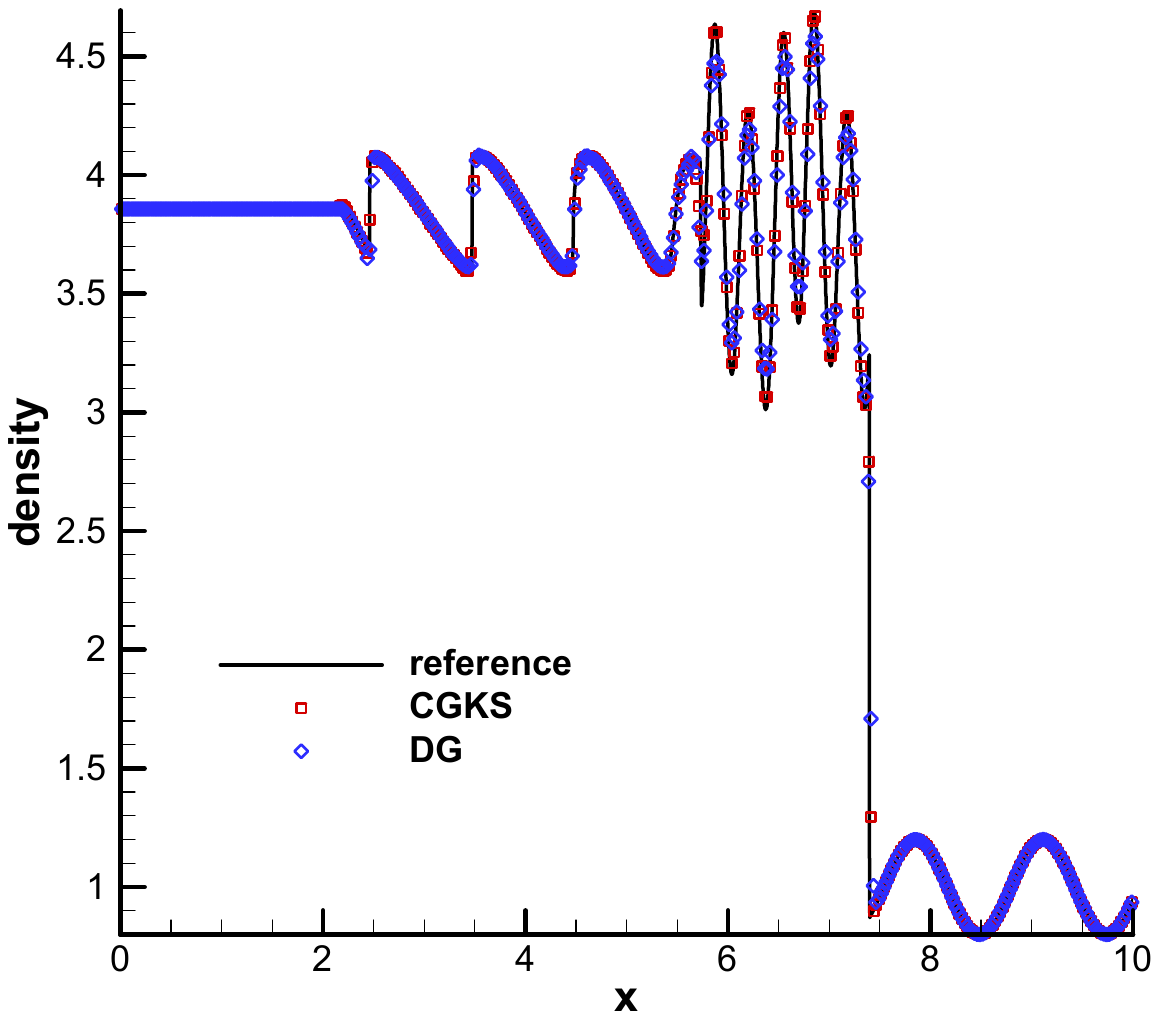}
		\includegraphics[width=0.32\textwidth, trim=20mm 0mm 30mm 10mm, clip]{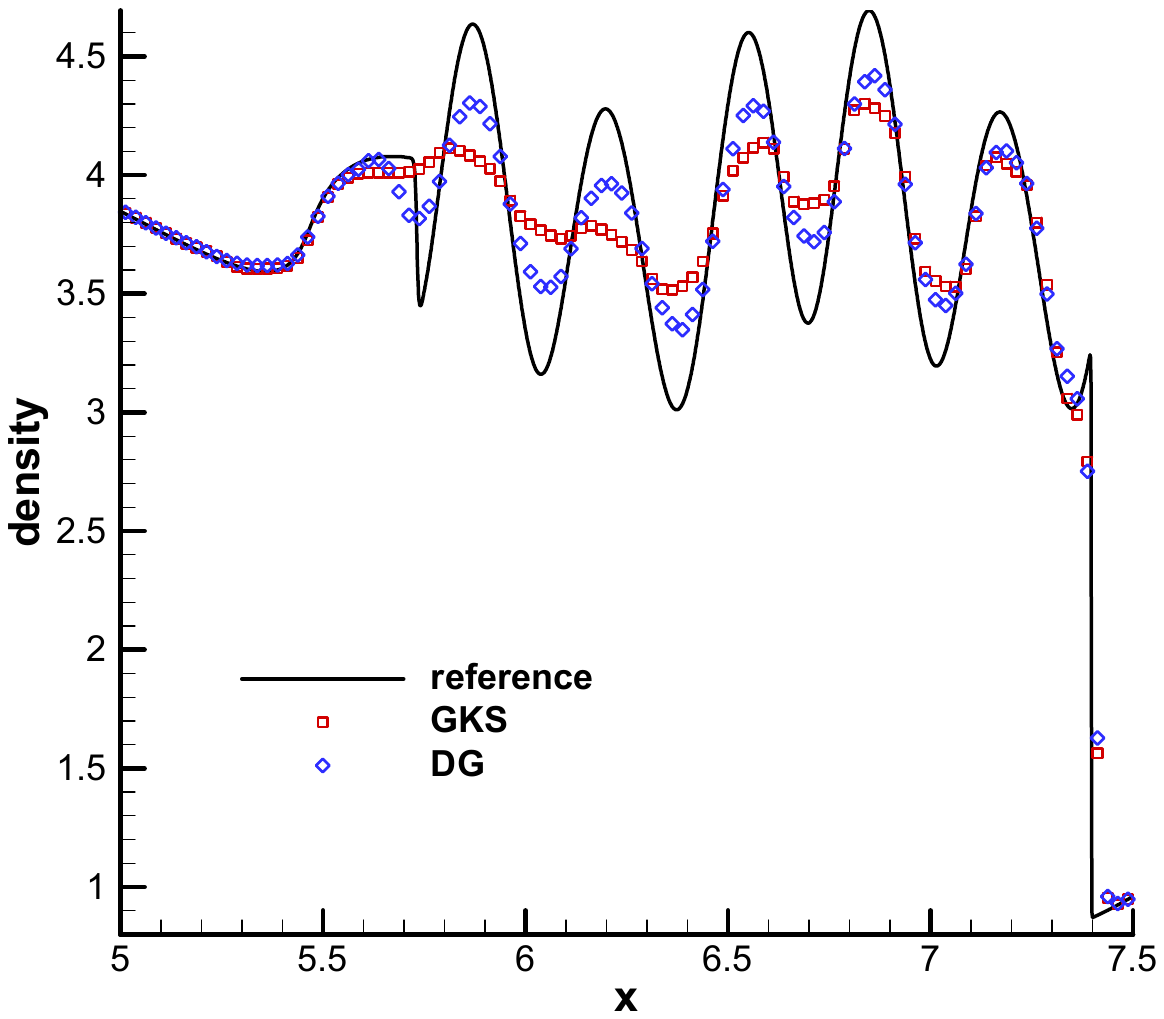}
		\includegraphics[width=0.32\textwidth, trim=20mm 0mm 30mm 10mm, clip]{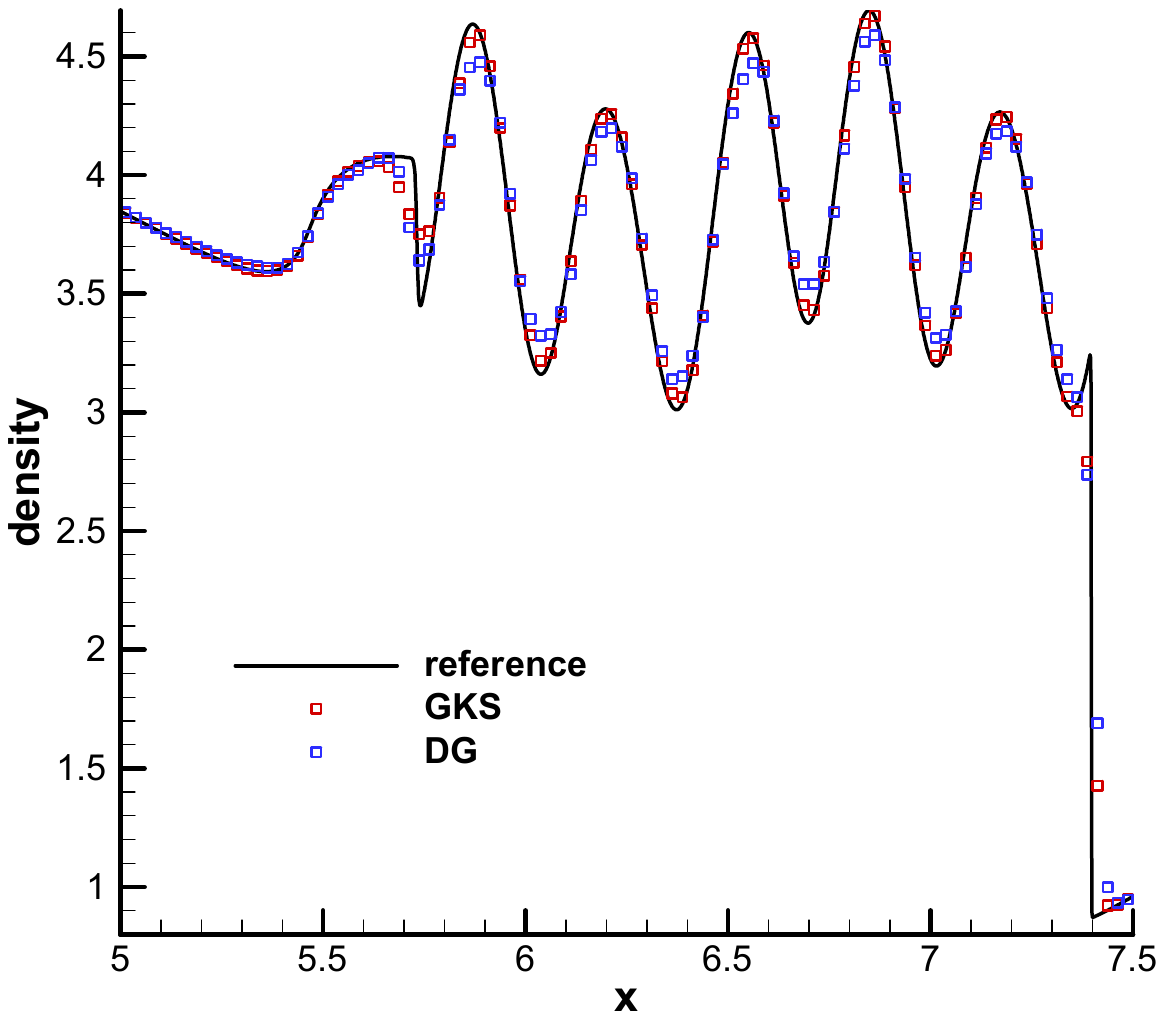}
		\includegraphics[width=0.32\textwidth, trim=20mm 0mm 30mm 10mm, clip]{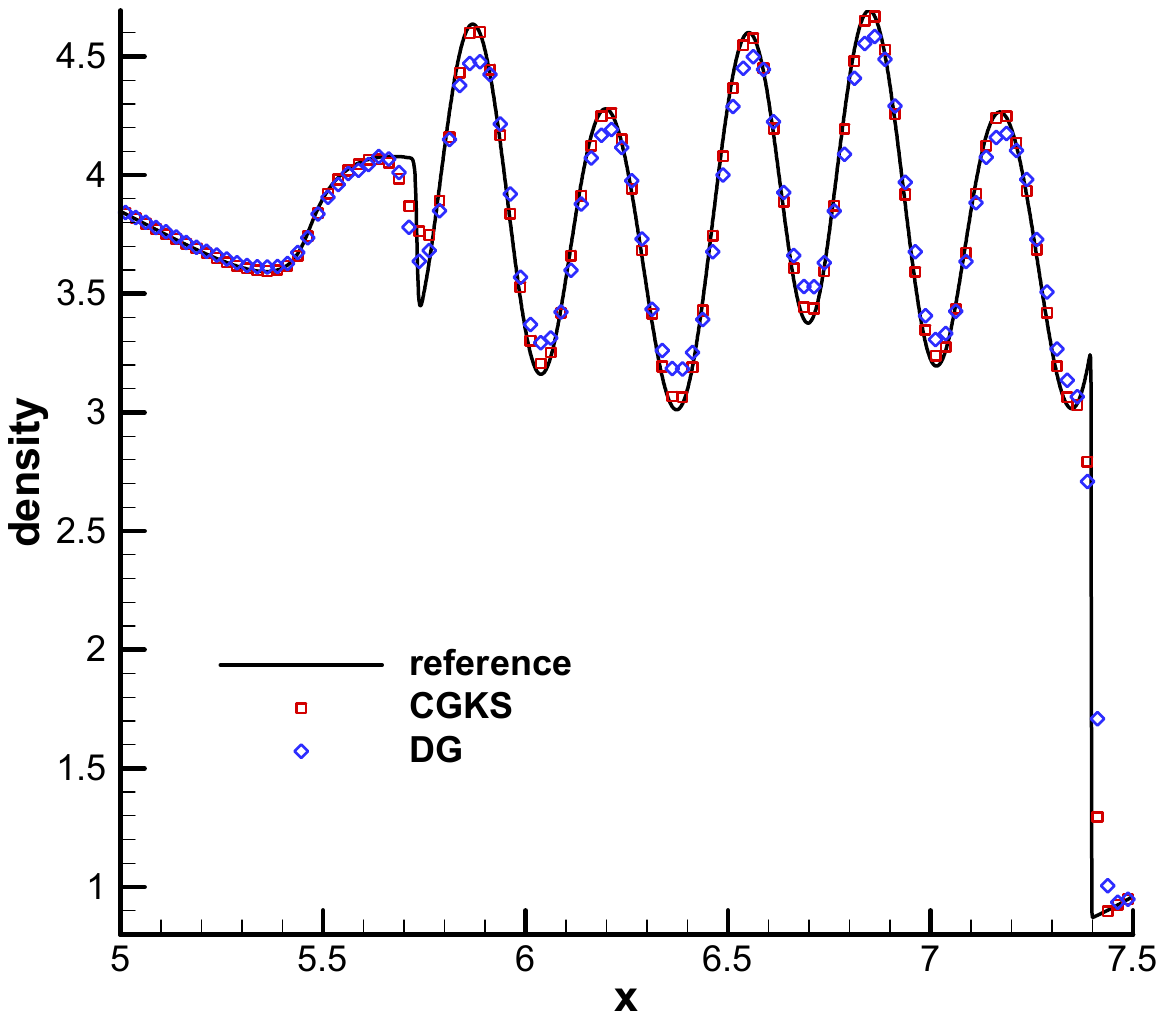}
		\caption{\label{shuosher} Shu-Osher problem: the density distributions and local enlargement at $t=1.8s$ with $400$ cells. Left to right: RKDG-$P^{2}$ with third-order GKS, RKDG-$P^{4}$ with fifth-order GKS, RKDG-$P^{4}$ with fifth-order CGKS.}
	\end{figure}

	\begin{figure}[!h]
		\centering
		\includegraphics[width=0.32\textwidth, trim=20mm 0mm 30mm 10mm, clip]{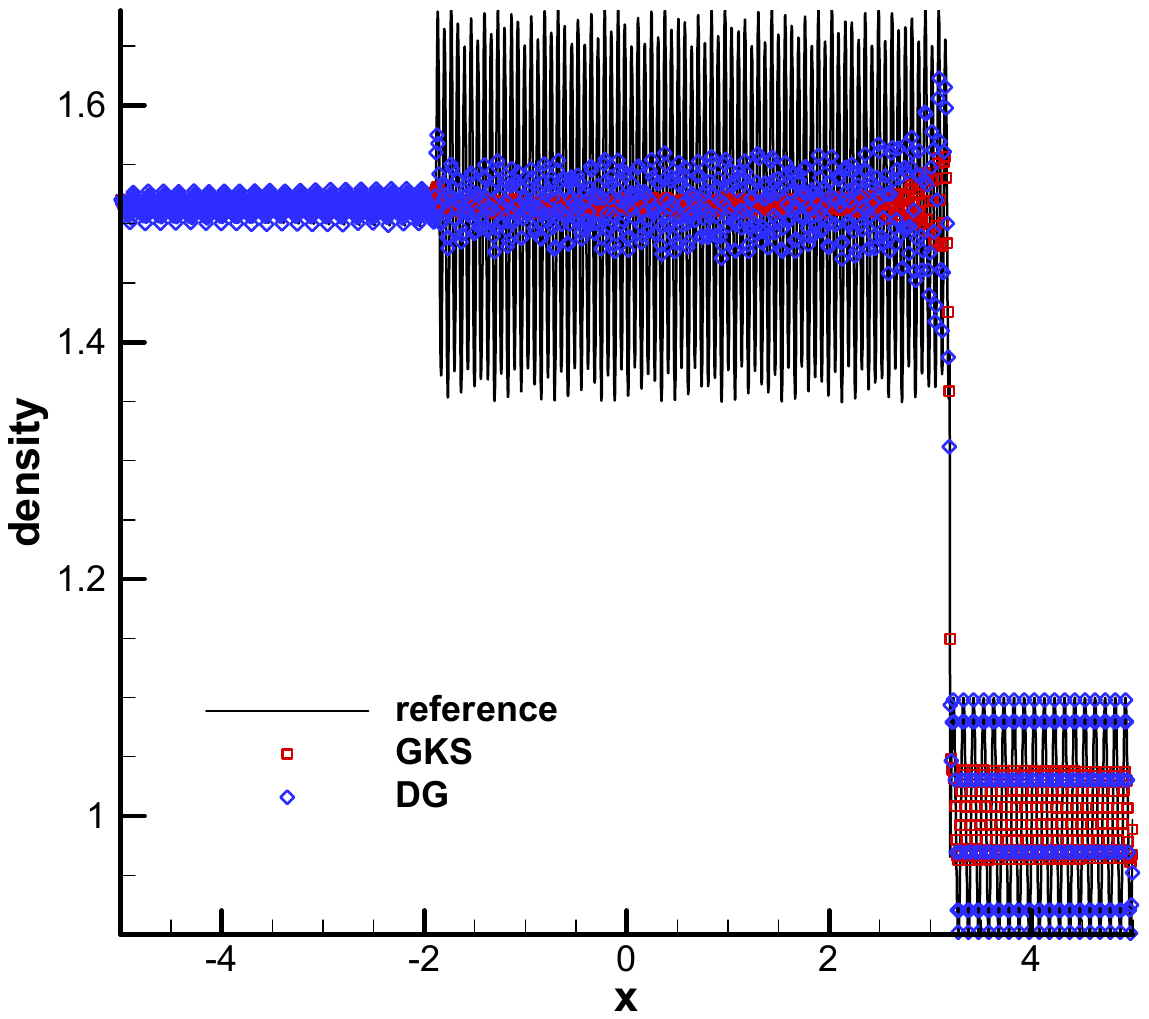}
		\includegraphics[width=0.32\textwidth, trim=20mm 0mm 30mm 10mm, clip]{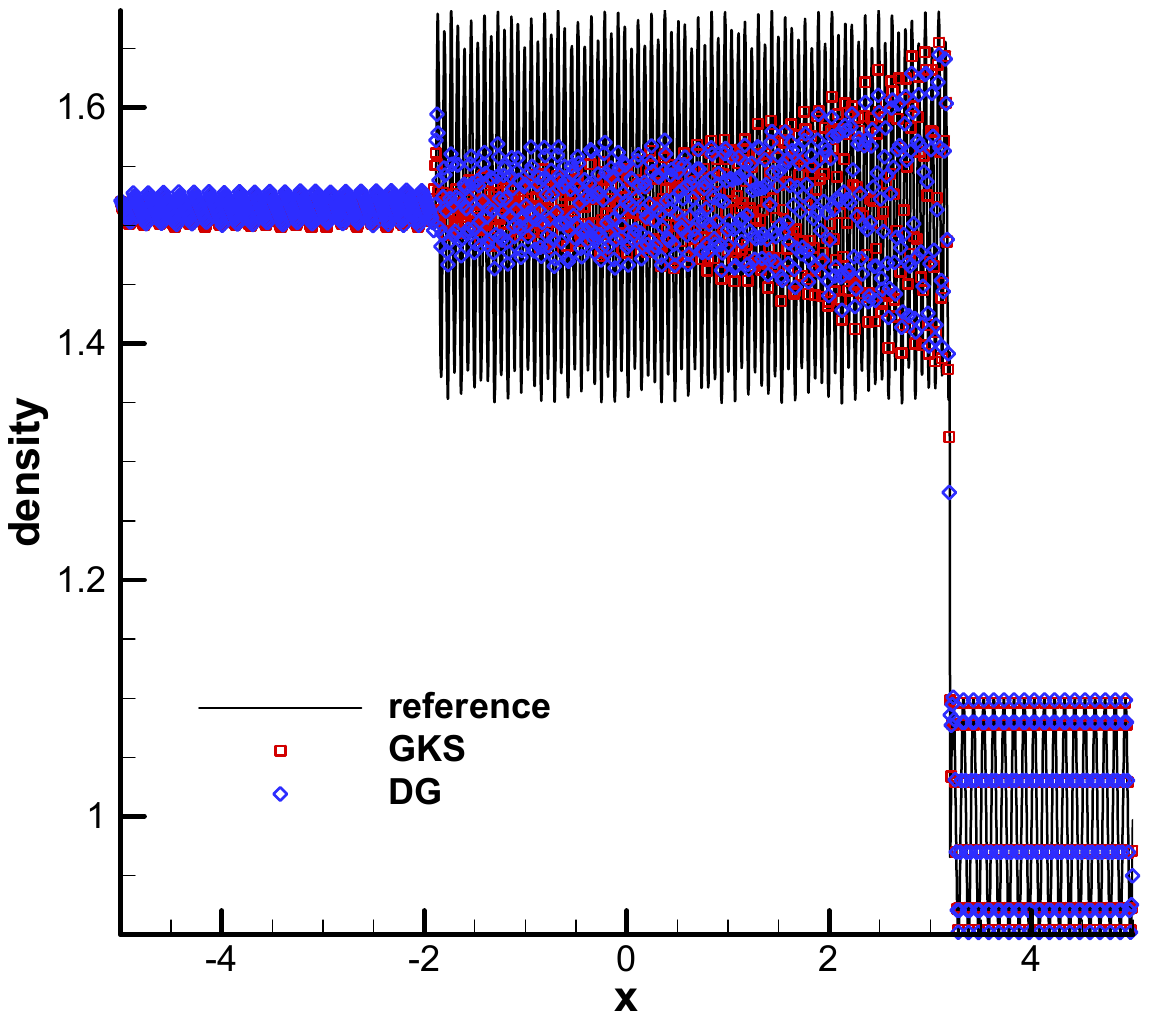}
		\includegraphics[width=0.32\textwidth, trim=20mm 0mm 30mm 10mm, clip]{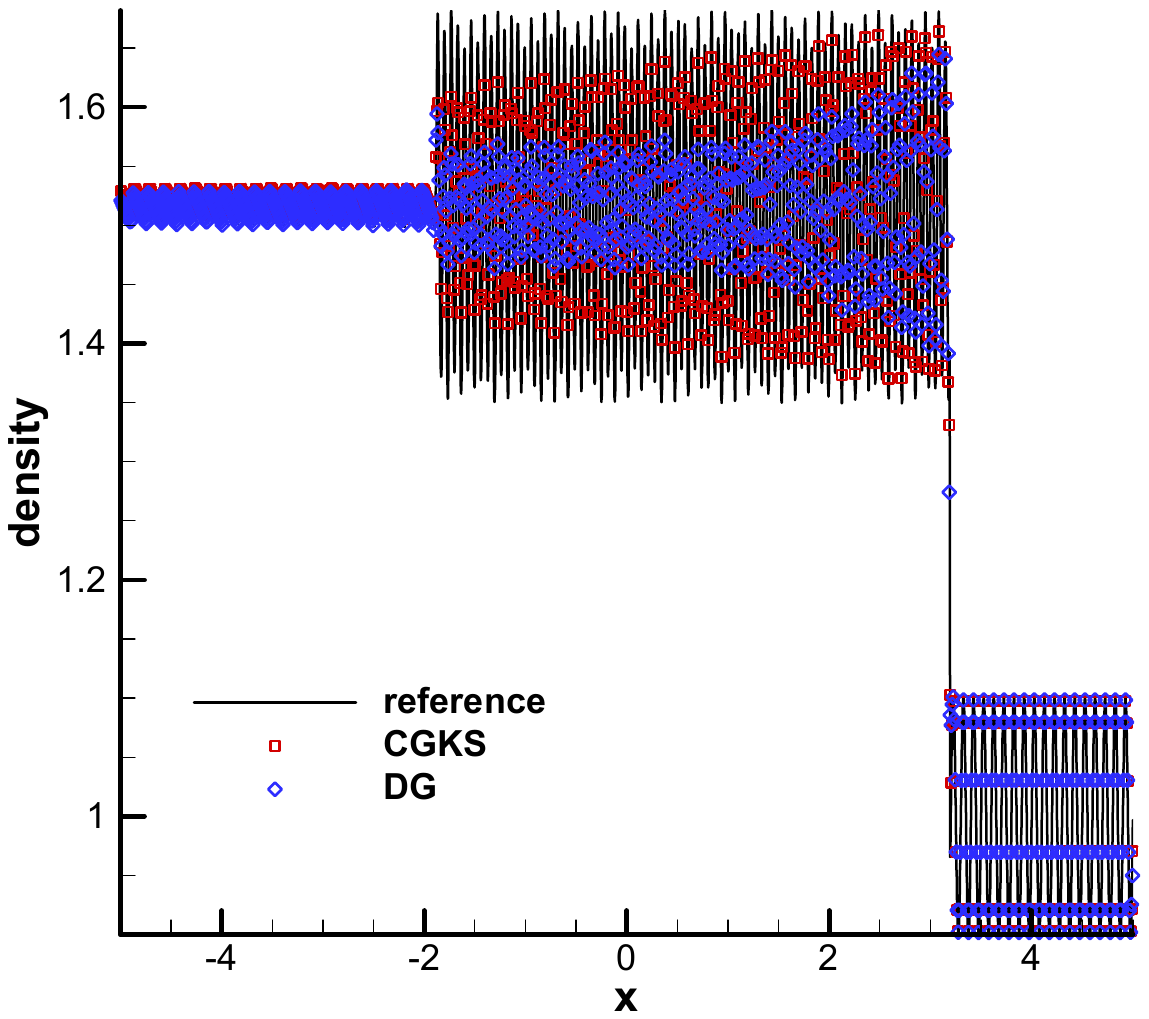}
		\includegraphics[width=0.32\textwidth, trim=20mm 0mm 30mm 10mm, clip]{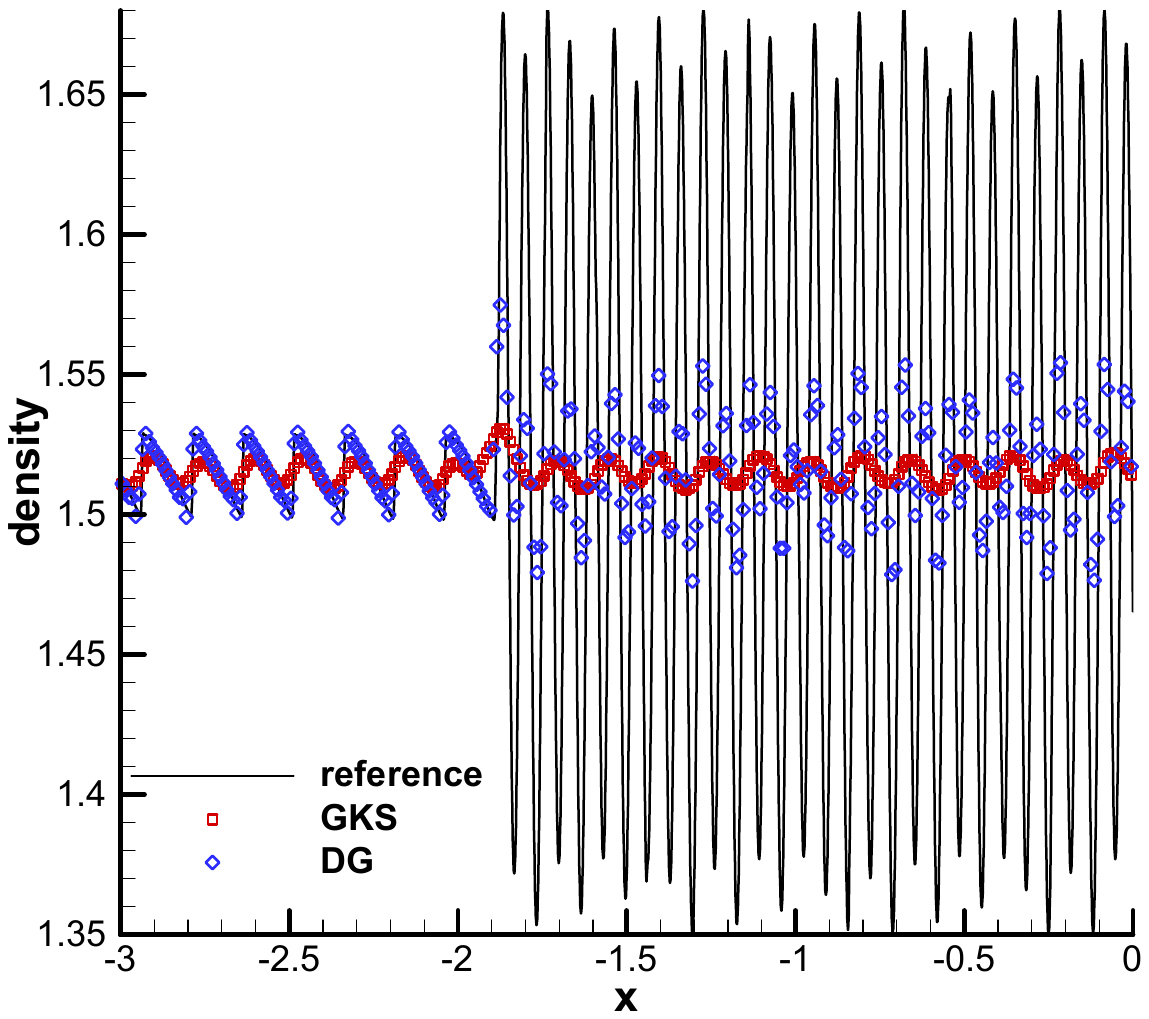}
		\includegraphics[width=0.32\textwidth, trim=20mm 0mm 30mm 10mm, clip]{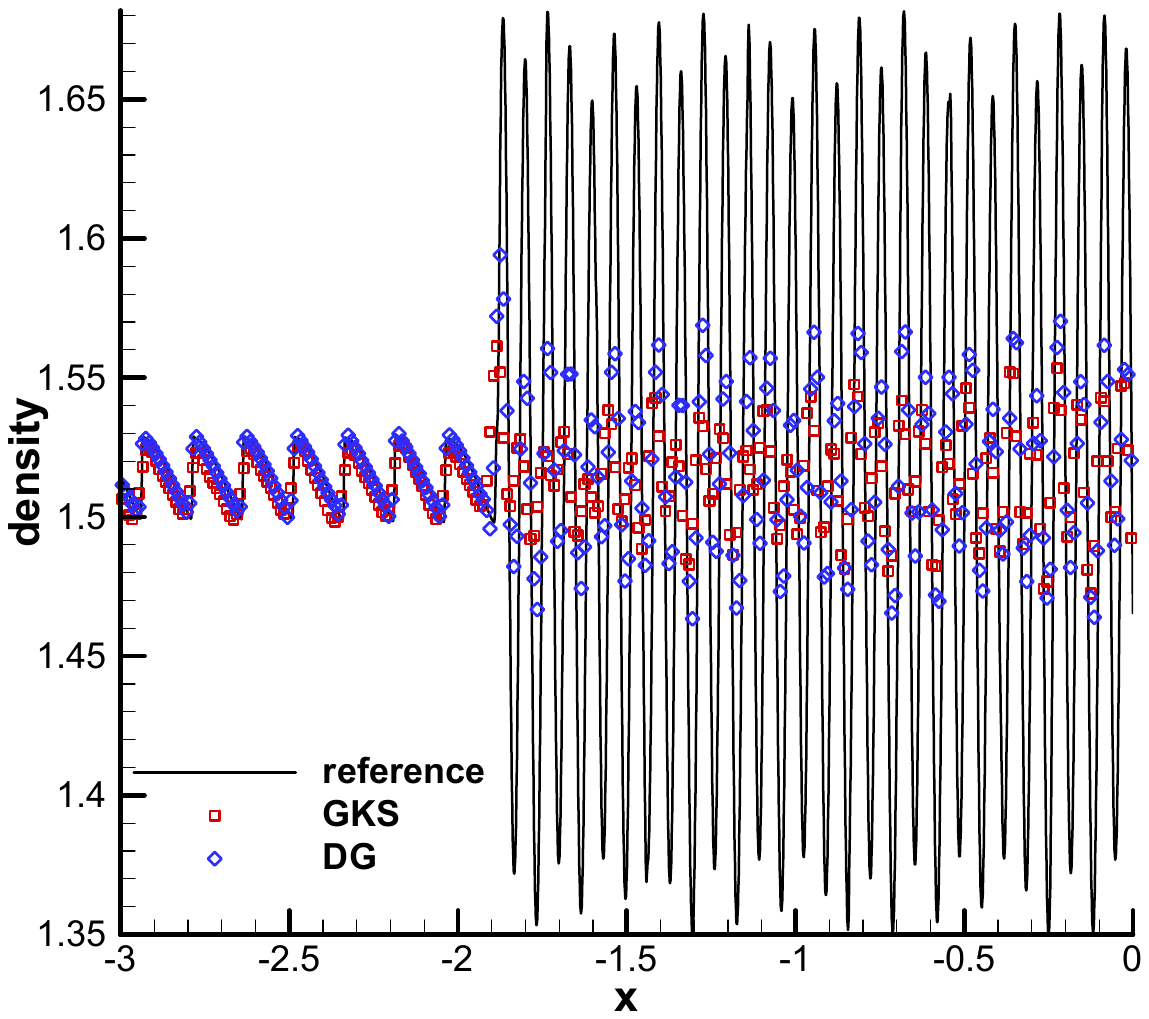}
		\includegraphics[width=0.32\textwidth, trim=20mm 0mm 30mm 10mm, clip]{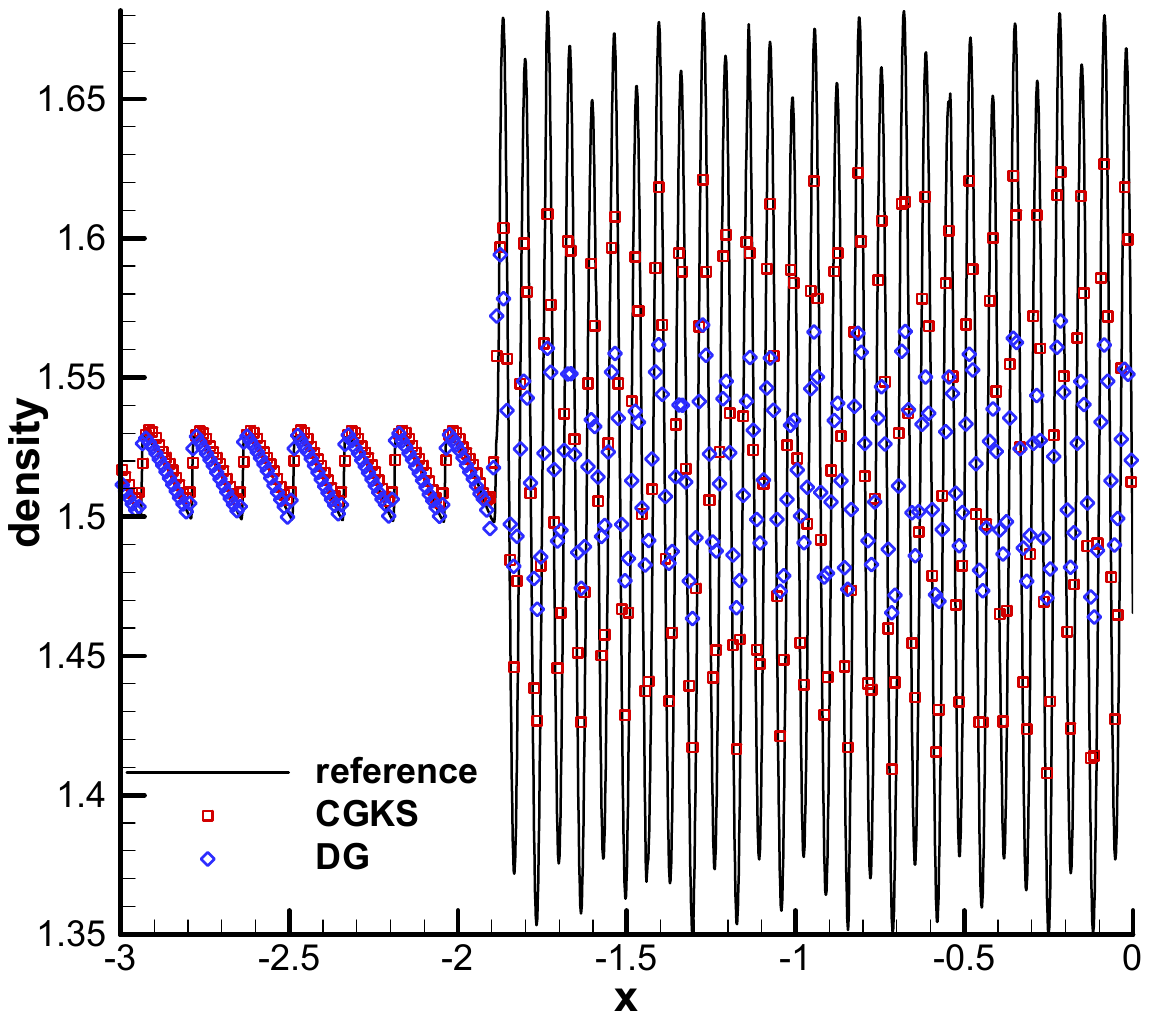}
		\caption{\label{toto} Titarev-Toro problem: the density distributions and local enlargement at $t=5s$ with $1000$ cells. Left to right: RKDG-$P^{2}$ with third-order GKS, RKDG-$P^{4}$ with fifth-order GKS, RKDG-$P^{4}$ with fifth-order CGKS.}
	\end{figure}

	\begin{figure}[!h]
		\centering
		\includegraphics[width=0.32\textwidth, trim=20mm 0mm 30mm 10mm, clip]{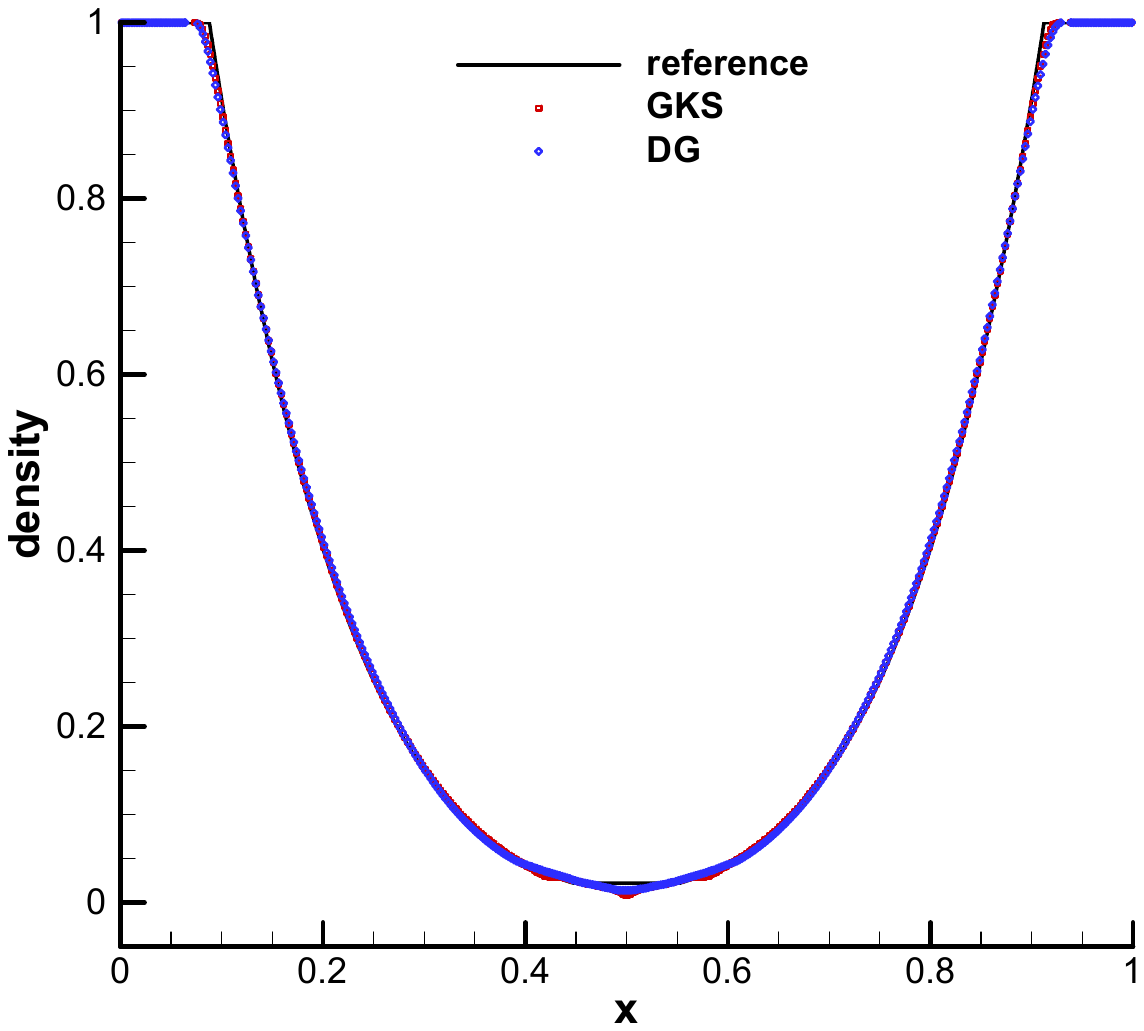}
		\includegraphics[width=0.32\textwidth, trim=20mm 0mm 30mm 10mm, clip]{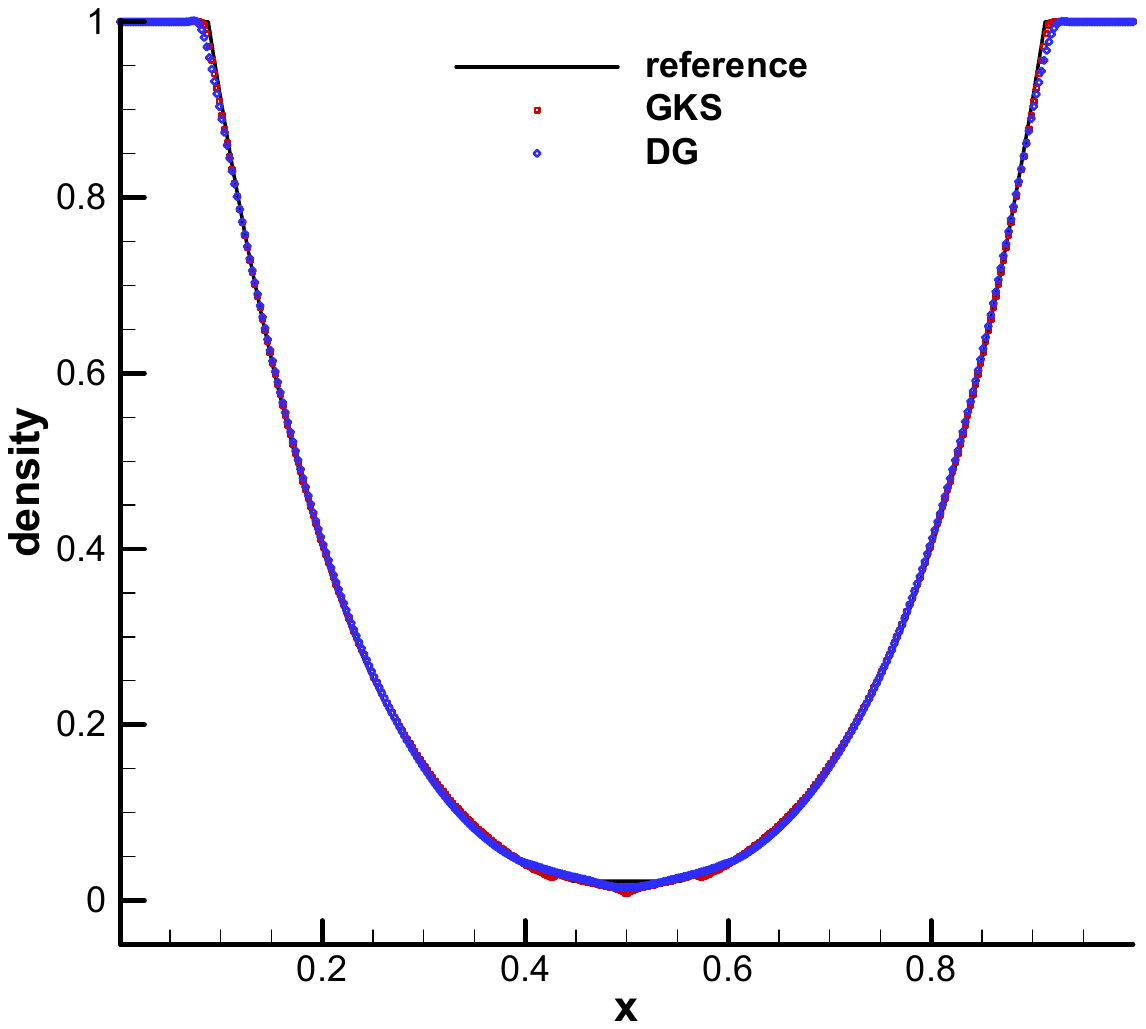}
		\includegraphics[width=0.32\textwidth, trim=20mm 0mm 30mm 10mm, clip]{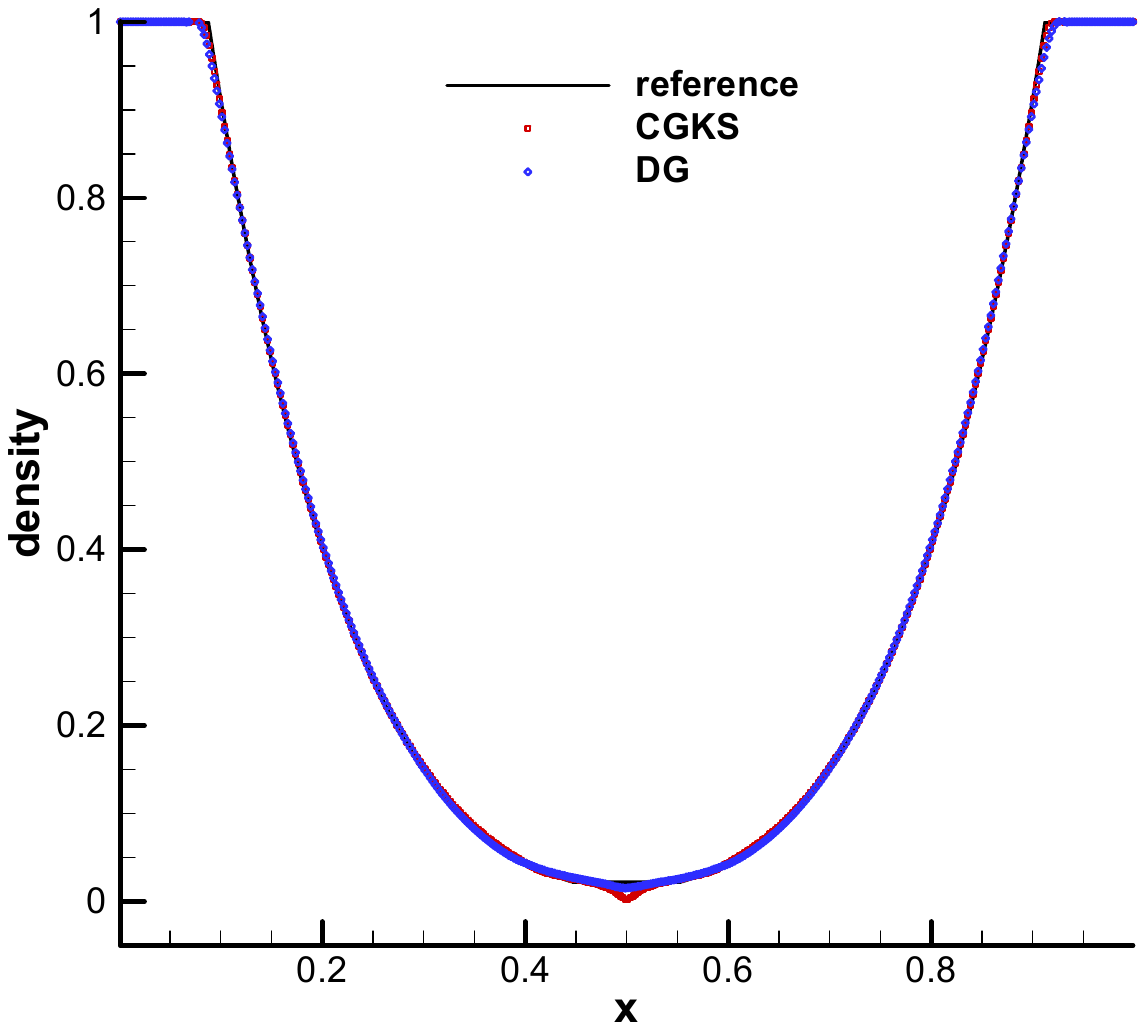}
		\caption{\label{rarewave} Double rarefaction waves: the density distributions and local enlargement at $t=0.15s$ with $400$ cells. Left to right: RKDG-$P^{2}$ with third-order GKS, RKDG-$P^{4}$ with fifth-order GKS, RKDG-$P^{4}$ with fifth-order CGKS.}
	\end{figure}

	\subsection{One dimensional problems with discontinuities}
	The first one is the Shu-Osher problem \cite{shu1988efficient}, which is used to test the performance of the numerical methods in a domain consisting of both shocks and complex smooth regions.
	The initial conditions are
	\begin{align*}
		(\rho,U,p) =\begin{cases}
			(3.857134, 2.629369, 10.33333), &  -5<x \leq -4,\\
			(1 + 0.2\sin (5x), 0, 1),  &  -4 <x<5,
		\end{cases}
	\end{align*}
	and the computational domain is (-5,5).
	The density distributions at $t=1.8s$ of the two methods against reference solution computed by the fifth-order WENO scheme with 10000 grid points are shown in Figure 4.
	
	As an extension of the Shu-Osher problem, the Titarev-Toro problem \cite{titarev2004finite} is tested as well, and
	the initial conditions in this case are the following
	\begin{align*}
		(\rho,U,p) =\begin{cases}
			(1.515695,0.523346,1.805),   & -5< x \leq -4.5,\\
			(1 + 0.1\sin (20\pi x), 0, 1),  &  -4.5 <x <5,
		\end{cases}
	\end{align*}
	and the computational domain is (-5,5). The
	non-reflecting boundary condition is imposed on left end, and the
	fixed wave profile is given on the right end in the two cases. The results of the two methods together with the reference solution is computed by fifth-order WENO scheme with 10000 grids are shown in Figure 5.
	As the results of the accuracy test, the performance of the RKDG method for simulating complex structure with shock waves and acoustic waves is between GKS  and CGKS.
	
	The last case for one-dimensional  problems is the double rarefaction wave problem \cite{linde1997robust}. The initial conditions are
	\begin{align}
		(\rho,U,p) =\begin{cases}
			(1, -2, 0.4), &  0<x \leq 0.5,\\
			(1, 2, 0.4),  &  0.5 <x<1,
		\end{cases}
	\end{align}
	and the computational domain is (0,1). The case is hard to simulate because of the near-vacuum region in the middle of domain. The results of the two methods together with exact solution computed by the exact Riemann solver \cite{toro1994restoration} are shown in Figure 6. What is worth mentioning is that in this case, the CFL number of the RKDG method is set to be 0.12 for RKDG-$P^{2}$ method and 0.04 for RKDG-$P^{4}$ method to avoid calculating a negative density. It is shown that the results simulated within GKS as well as CGKS  have observable numerical  oscillations in the mid region, where the density is at its minimum. However, the results calculated by RKDG have smaller numerical oscillations as the ones calculated by GKS (CGKS).

	\subsection{Two dimensinal accuracy test}
	In this subsection, two dimensional advection of density perturbation is tested, and the initial
	condition is given as follows
	\begin{align*}
		\rho(x,y)=1+0.2\sin(\pi (x+y)),\ \  U(x,y)=1,\ \
		\ V(x,y)=1 ,\ \ p(x,y)=1, x,y\in[-1,1].
	\end{align*}
	With the periodic boundary condition, and the analytic
	solution is
	\begin{align*}
		\rho(x,y,t)=1+0.2\sin(\pi (x+y-t)),\ \  U(x,y,t)=1,\ \
		\ V(x,y,t)=1 ,\ \ p(x,y,t)=1.
	\end{align*}
	The simulation time is $t=2s$ in this case. And fourth-order Runge-Kutta method \cite{butcher1996history} is adopted for RKDG-$P^{4}$   to guarantee that the spatial error dominates. The $L^{1}$ and $L^{2}$ errors together with their orders are shown in Table 5-6 for the fifth-order RKDG and GKS(CGKS). Also, the comparison of the efficiency is presented in Figure 7.
	
	It is shown that the RKDG-$P^{4}$ method does perform better with extremely small errors compared to the fifth-order GKS(CGKS) in the same uniform meshes. As for the efficiency, RKDG-$P^{4}$ and the CGKS have similar performance, which is much better than the GKS's.
	\begin{figure}[!h]
		\centering
		\includegraphics[width=0.49\textwidth, trim=20mm 0mm 30mm 10mm, clip]{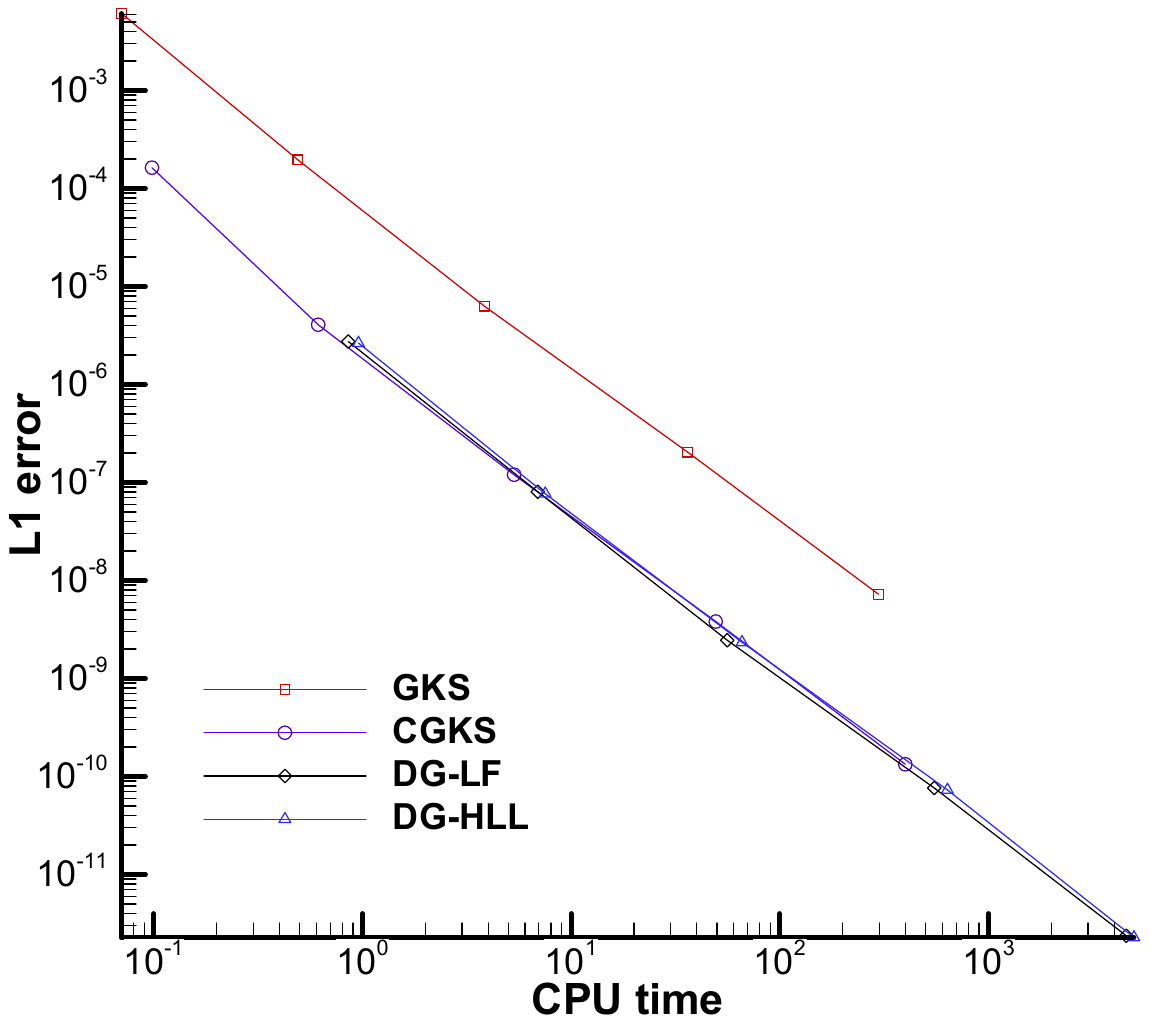}
		\includegraphics[width=0.49\textwidth, trim=20mm 0mm 30mm 10mm, clip]{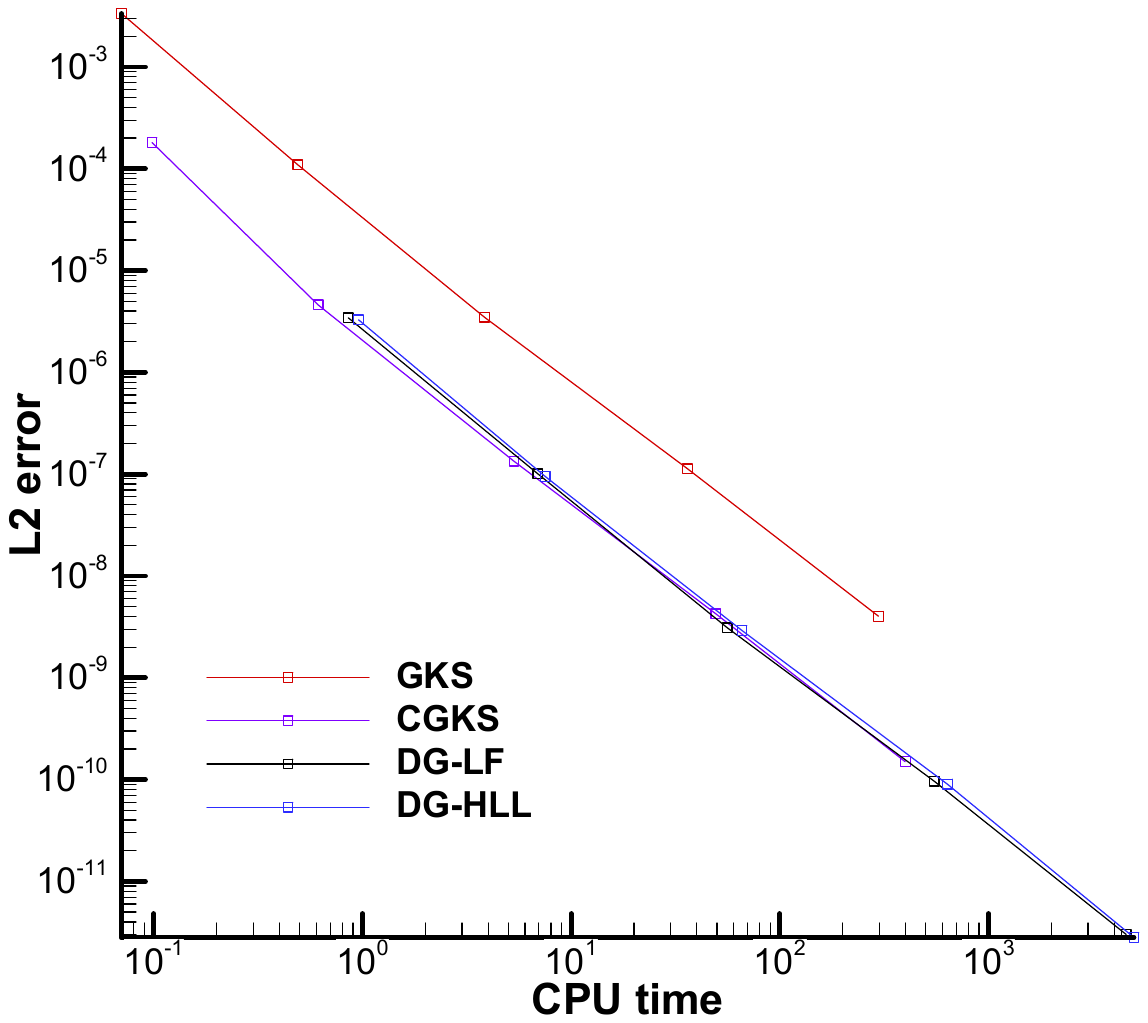}
		\caption{\label{accuracy-5-td} Efficiency of simulation for two dimensional smooth regions:comparisons between RKDG-$P^4$ method and fifth-order GKS/CGKS method. The left one is based on $L^{1}$ error and the right one is based on $L^{2}$ error}.
	\end{figure}
	
	\begin{table}[]
		\centering
		\begin{tabular}{|l|llll|llll|}
			\hline
			Mesh & \multicolumn{4}{c|}{RKDG-$P^{4}$ with LF flux}                                                          & \multicolumn{4}{c|}{RKDG-$P^{4}$ with HLL flux}                                                         \\ \hline
			& \multicolumn{1}{l|}{$L^{1} $ error} & \multicolumn{1}{l|}{order} & \multicolumn{1}{l|}{$L^{2} $ error} & order & \multicolumn{1}{l|}{$L^{1} $ error} & \multicolumn{1}{l|}{order} & \multicolumn{1}{l|}{$L^{2} $ error} & order \\ \hline
			10   & \multicolumn{1}{l|}{2.75E-06} & \multicolumn{1}{l|}{}      & \multicolumn{1}{l|}{3.43E-06} &       & \multicolumn{1}{l|}{2.65E-06} & \multicolumn{1}{l|}{}      & \multicolumn{1}{l|}{3.28E-06} &       \\ \hline
			20   & \multicolumn{1}{l|}{8.06E-08} & \multicolumn{1}{l|}{5.09}  & \multicolumn{1}{l|}{1.01E-07} & 5.08  & \multicolumn{1}{l|}{7.69E-08} & \multicolumn{1}{l|}{5.11}  & \multicolumn{1}{l|}{9.51E-08} & 5.11  \\ \hline
			40   & \multicolumn{1}{l|}{2.45E-09} & \multicolumn{1}{l|}{5.04}  & \multicolumn{1}{l|}{3.10E-09} & 5.03  & \multicolumn{1}{l|}{2.34E-09} & \multicolumn{1}{l|}{5.04}  & \multicolumn{1}{l|}{2.91E-09} & 5.03  \\ \hline
			80   & \multicolumn{1}{l|}{7.61E-11} & \multicolumn{1}{l|}{5.01}  & \multicolumn{1}{l|}{9.64E-11} & 5.01  & \multicolumn{1}{l|}{7.27E-11} & \multicolumn{1}{l|}{5.01}  & \multicolumn{1}{l|}{9.02E-11} & 5.01  \\ \hline
			160  & \multicolumn{1}{l|}{2.37E-12} & \multicolumn{1}{l|}{5.00}  & \multicolumn{1}{l|}{3.01E-12} & 5.00  & \multicolumn{1}{l|}{2.27E-12} & \multicolumn{1}{l|}{5.00}  & \multicolumn{1}{l|}{2.81E-12} & 5.00  \\ \hline
		\end{tabular}
		\caption{Two dimensional advection of density perturbation: accuracy test for RKDG-$P^{4}$ with LF flux and HLL flux without limiters}
		\label{tab:my-table}
	\end{table}
	\begin{table}[]
		\centering
		\begin{tabular}{|l|llll|llll|}
			\hline
			Mesh & \multicolumn{4}{c|}{Fifth-order GKS}                                                               & \multicolumn{4}{c|}{Fifth-order CGKS}                                                              \\ \hline
			& \multicolumn{1}{l|}{$L^{1} $ error} & \multicolumn{1}{l|}{order} & \multicolumn{1}{l|}{$L^{2} $ error} & order & \multicolumn{1}{l|}{$L^{1} $ error} & \multicolumn{1}{l|}{order} & \multicolumn{1}{l|}{$L^{2} $ error} & order \\ \hline
			10   & \multicolumn{1}{l|}{6.11E-03} & \multicolumn{1}{l|}{}      & \multicolumn{1}{l|}{3.35E-03} &       & \multicolumn{1}{l|}{1.63E-04} & \multicolumn{1}{l|}{}      & \multicolumn{1}{l|}{1.81E-04} &       \\ \hline
			20   & \multicolumn{1}{l|}{1.97E-04} & \multicolumn{1}{l|}{4.96}  & \multicolumn{1}{l|}{1.09E-04} & 4.94  & \multicolumn{1}{l|}{4.07E-06} & \multicolumn{1}{l|}{5.33}  & \multicolumn{1}{l|}{4.60E-06} & 5.29  \\ \hline
			40   & \multicolumn{1}{l|}{6.28E-06} & \multicolumn{1}{l|}{4.97}  & \multicolumn{1}{l|}{3.48E-06} & 4.97  & \multicolumn{1}{l|}{1.20E-07} & \multicolumn{1}{l|}{5.09}  & \multicolumn{1}{l|}{1.34E-07} & 5.10  \\ \hline
			80   & \multicolumn{1}{l|}{2.04E-07} & \multicolumn{1}{l|}{4.94}  & \multicolumn{1}{l|}{1.13E-07} & 4.94  & \multicolumn{1}{l|}{3.80E-09} & \multicolumn{1}{l|}{4.98}  & \multicolumn{1}{l|}{4.26E-09} & 4.98  \\ \hline
			160  & \multicolumn{1}{l|}{7.21E-09} & \multicolumn{1}{l|}{4.82}  & \multicolumn{1}{l|}{4.01E-09} & 4.82  & \multicolumn{1}{l|}{1.33E-10} & \multicolumn{1}{l|}{4.83}  & \multicolumn{1}{l|}{1.50E-10} & 4.82  \\ \hline
		\end{tabular}
		\caption{Two dimensional advection of density perturbation: accuracy test for fifth-order GKS and CGKS method with smooth reconstruction}
		
	\end{table}

	\subsection{Double Mach reflection problem}
	This problem was extensively studied by Woodward and Colella
	\cite{woodward1984numerical} for the inviscid flow. The computational
	domain is $[0,4]\times[0,1]$, and a solid wall lies at the bottom of
	the computational domain starting from $x =1/6$. Initially, a
	right-moving Mach $10$ shock is positioned at $(x,y)=(1/6, 0)$, and
	makes a $60^\circ$ angle with the x-axis. The initial pre-shock and
	post-shock conditions are
	\begin{align*}
		(\rho, U, V, p)&=(8, 4.125\sqrt{3}, -4.125,
		116.5),\\
		(\rho, U, V, p)&=(1.4, 0, 0, 1).
	\end{align*}
	The reflecting boundary condition is used at the wall. While for the
	rest of the bottom boundary, the exact post-shock condition is imposed.
	At the top boundary, the flow variables are set to follow  the
	motion of the Mach $10$ shock. In this case, we tend to compare the resolution of the two methods under the premise of almost the same computational load. A baseline uniform mesh with $800\times200$ grid points is used for the two methods, and a refined one ($1600\times400$) is used in the CGKS method to make the comparison fair. There are $15$ degrees of freedom contained in a single cell in the RKDG-$P{4}$ method, but only one and three degrees of freedom for each variable in the GKS and CGKS, respectively. Thus, the refined one has similar degrees of freedom($1,920,000$ totally) with the RKDG-$P{4}$ method($2,400,000$ totally). As is shown in Figure 7 and Figure 8, the CGKS(refined one) has a similar resolution with the RKDG-$P^{4}$ method. However, discrete local contours are contained in the RKDG method even with the application of the trouble cell indicators. The spurious oscillations are much smaller in the CGKS, which validates its better fidelity on simulating high-speed compressible flow.
	
	In this case, 16 CPU cores are used parallelly and the time consumption for the first 0.02s simulation time is shown in Table \ref{CPU-double-mach}. Under approximately the same degrees of freedom, RKDG and CGKS seem to have similar resolutions, while the computational overhead of CGKS is smaller.
	
	\begin{figure}
		\centering
		\includegraphics[width=1\textwidth, trim=0mm 110mm 0mm 114.5mm, clip]{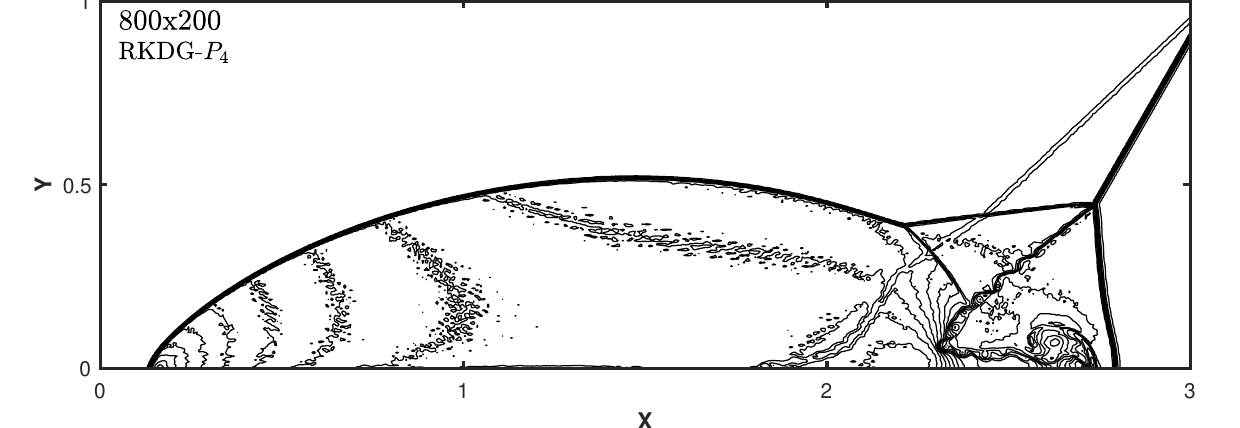}
		\includegraphics[width=1\textwidth, trim=0mm 110mm 0mm 114.5mm, clip]{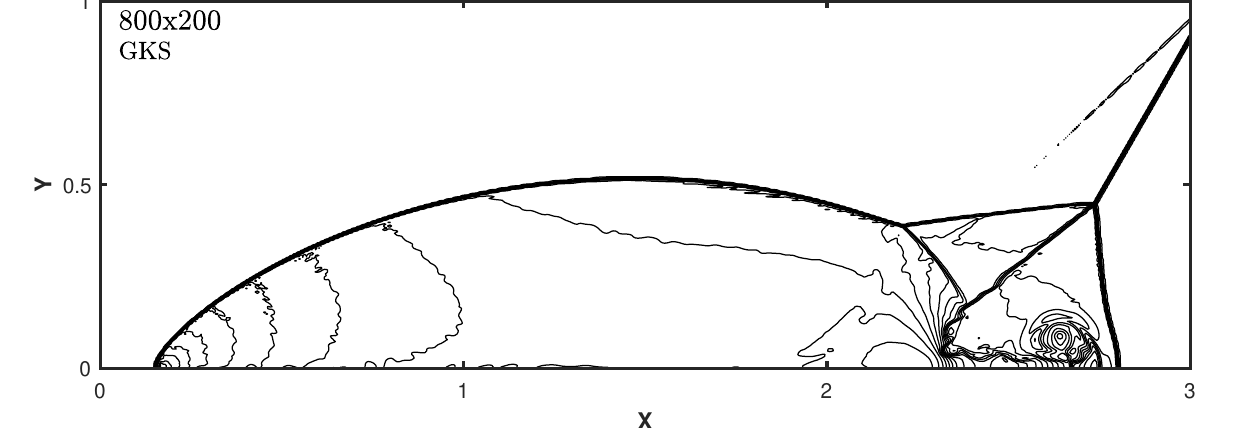}
		\includegraphics[width=1\textwidth, trim=0mm 110mm 0mm 114.5mm, clip]{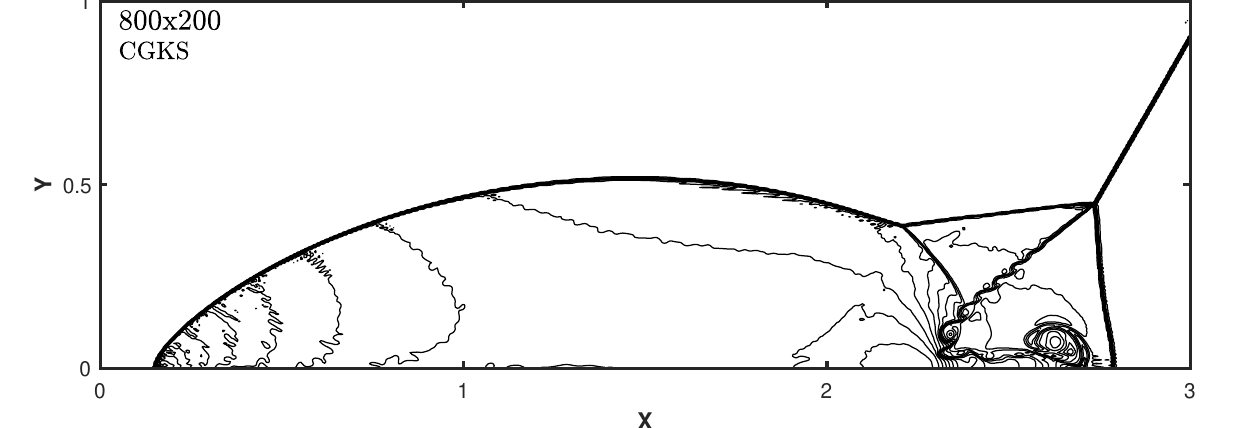}
		\includegraphics[width=1\textwidth, trim=0mm 110mm 0mm 114.5mm, clip]{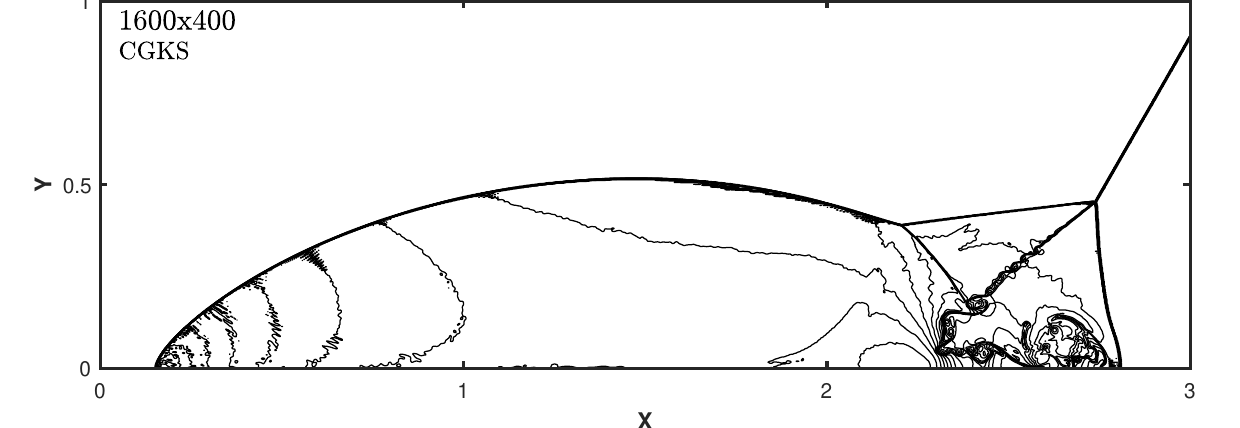}

		\caption{\label{DM} Double Mach problems: the density distributions at $t=0.2s$. 30 equally spaced density contours from 1.5 to 21.5. From top to bottom: RKDG-$P^{4}$, fifth-order GKS, fifth-order CGKS, fifth-order CGKS(refined).}
	\end{figure}
	
	\begin{figure}
		\centering
		\includegraphics[width=0.39\textwidth, trim=50mm 110mm 50mm 114.5mm, clip]{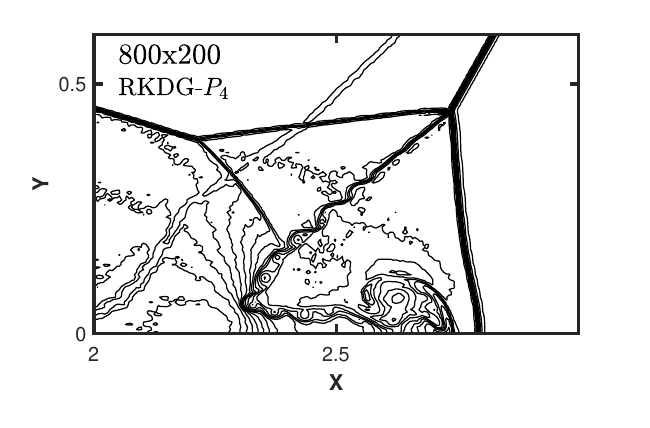}
		\includegraphics[width=0.39\textwidth, trim=50mm 110mm 50mm 114.5mm, clip]{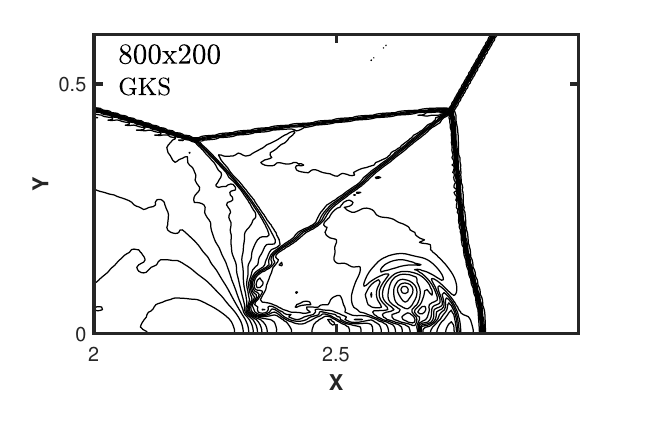}
		\includegraphics[width=0.39\textwidth, trim=50mm 110mm 50mm 114.5mm, clip]{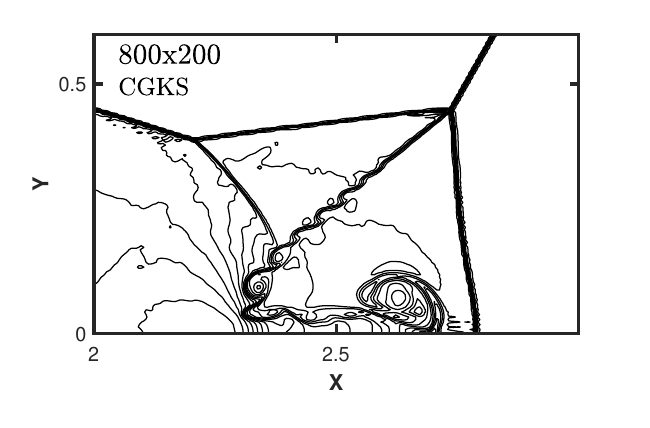}
		\includegraphics[width=0.39\textwidth, trim=50mm 110mm 50mm 114.5mm, clip]{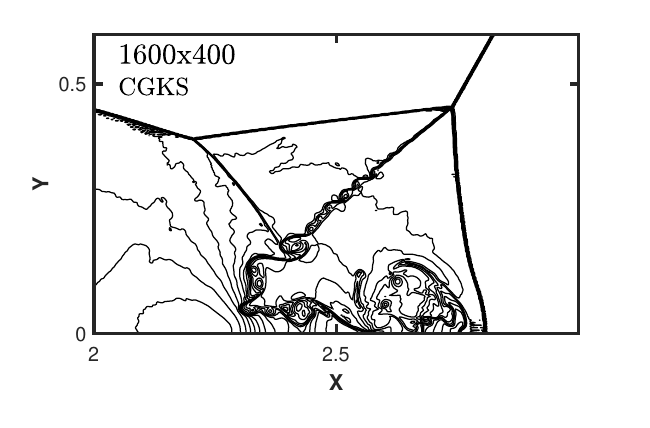}

		\caption{\label{DM_zoom} Double Mach problems: zoom-in pictures around the Mach stem. 30 equally spaced density contours from 1.5 to 21.5.From left to right and top to bottom: RKDG-$P^{4}$, fifth-order GKS, fifth-order CGKS, fifth-order CGKS(refined).}

	\end{figure}
	
		\begin{table}[]
		\centering
		\begin{tabular}{|l|l|l|l|}
			\hline
			Scheme & RKDG-$P^{4}$    & CGKS     & CGKS        \\ \hline
			Mesh Number & 800$\times$200 & 800$\times$200 & 1600$\times$400
			\\ \hline
			CPU time & {749.179s} & {75.6s} & {588.577s} \\ \hline
		\end{tabular}
		\caption{Double Mach problems: the CPU time consumption for the first 0.02s simulation time}
		\label{CPU-double-mach}
	\end{table}
	
	\section{Conclusion}

In this paper, we systematically evaluate the performance of the RKDG method and the GKS in terms of accuracy, efficiency, and resolution. We demonstrate significant improvements in the efficiency of the GKS through the use of compact spatial reconstruction and a simplified third-order gas distribution method. Notably, GKS exhibits higher efficiency than the RKDG method in smooth flow regions. Importantly, GKS is capable of handling viscous flows without additional computational demands, whereas implementing viscous terms in the RKDG method requires further development, indicating a potentially higher efficiency of the CGKS for real-world viscous flow simulations. 
During the implementation of these methods, we encountered challenges in maintaining the stability of the RKDG method when using multi-resolution WENO limiters. The coefficients $\gamma_{\ell-1, \ell}$ and $\gamma_{\ell, \ell}$ 
for the limiters proved sensitive and required case-specific adjustments to achieve both robustness and accuracy in capturing discontinuous solutions. Furthermore, the multi-resolution scheme tends to incorrectly identify excessive extreme points in smooth regions, thereby attenuating extreme values. Our comparative analysis highlights the robustness of the CGKS and the high accuracy of the RKDG method in multi-dimensional scenarios. Thus, integrating the strengths of GKS's relaxed CFL conditions and DG's high accuracy, a hybrid scheme combining DG and CGKS with $h-p$ 
refinement emerges as a promising approach for effectively simulating flows in both smooth and discontinuous regions. 
This development is part of our ongoing research efforts.
	
\section*{Acknowledgment}
	The current research is supported by National Science Foundation of China (12302378, 92371201, 92371107,12172316),
	Hong Kong Research Grant Council (16208021,16301222).

	\bibliographystyle{plain}
	\bibliography{bibt.bib}
	
\end{document}